\documentclass[a4paper]{article}
\def\({\left(}
\def\){\right)}
\usepackage{amsthm,amsmath,amssymb,color}
\newtheorem{teor}{Theorem}[section]
\newtheorem{defin}[teor]{Definition}
\newtheorem{lemm}[teor]{Lemma}
\newtheorem{prop}[teor]{Proposition}
\newtheorem{coro}[teor]{Corollary}
\newtheorem{prob}[teor]{Problem}
\theoremstyle{definition}
\newtheorem{defi}[teor]{Definition}
\newtheorem{osse}[teor]{Remark}

\newcommand{\bele}{\begin{lemm}\begin{sl}}
\newcommand{\enle}{\end{sl}\end{lemm}}
\newcommand{\bedef}{\begin{defi}\begin{sl}}
\newcommand{\eddef}{\end{sl}\end{defi}}
\newcommand{\bete}{\begin{teor}\begin{sl}}
\newcommand{\ente}{\end{sl}\end{teor}}
\newcommand{\beos}{\begin{osse}\begin{rm}}
\newcommand{\eddos}{\end{rm}\end{osse}}
\newcommand{\bepr}{\begin{prop}\begin{sl}}
\newcommand{\empr}{\end{sl}\end{prop}}
\newcommand{\bepro}{\begin{prob}\begin{rm}}
\newcommand{\empro}{\end{rm}\end{prob}}
\newcommand{\bede}{\begin{defin}\begin{sl}}
\newcommand{\edde}{\end{sl}\end{defin}}
\newcommand{\beco}{\begin{coro}\begin{sl}}
\newcommand{\enco}{\end{sl}\end{coro}}
\newcommand{\disp}{\displaystyle}

\newcommand{\thspace}{\hspace{3mm}}

\newcommand{\RR}{\mathbb{R}}

\newcommand{\R}{\mathbb R}
\newcommand{\eb}{\varepsilon}
\newcommand{\Dt}{\partial_t}

\newcommand{\beeq}[1]{\begin{equation}\label{#1}}
\newcommand{\eddeq}{\end{equation}}

\newcommand{\beeqa}[1]{\begin{eqnarray}\label{#1}}
\newcommand{\eddeqa}{\end{eqnarray}}

\newcommand{\beal}[1]{\begin{align}\label{#1}}
\newcommand{\eddal}{\end{align}}

\newcommand{\bespl}[1]{\begin{split}\label{#1}}
\newcommand{\edspl}{\end{split}}

\newcommand{\bega}[1]{\begin{gather}\label{#1}}
\newcommand{\edga}{\end{gather}}

\newcommand{\beeqax}{\begin{eqnarray*}}
\newcommand{\eddeqax}{\end{eqnarray*}}

\def\qed{\ifmmode 
  \else \leavevmode\unskip\penalty9999 \hbox{}\nobreak\hfill
  \fi
  \quad\hbox{\hskip.5em\vrule width.4em height.6em depth.05em\hskip.1em}}
\def\endproofsym{\qed}
\renewenvironment{proof}[1][Proof]{\trivlist\item[\hskip\labelsep{\hskip0pt
    {\normalfont\scshape#1.}\hskip .321429\parindent}]\ignorespaces}
{\endproofsym\endtrivlist}
\def\endnobox{\def\endproofsym{}\end{proof}\def\endproofsym{\qed}}

\newcommand{\no}{\nonumber}

\newcommand{\beeqao}{\begin{eqnarray}\no}
\newcommand{\bealo}{\begin{align}\no}
\newcommand{\besplo}{\begin{split}\no}
\newcommand{\begao}{\begin{gather}\no}

\newcommand{\eps}{\varepsilon}

\newcommand{\duav}[1]{\langle{#1}\rangle}

\newcommand{\dtt}{_{tt}}

\newcommand{\perogni}{\forall\,}
\newcommand{\esiste}{\exists\,}

\newcommand{\io}{\int_\Omega}

\newcommand{\epsi}{\varepsilon}
\newcommand{\ee}{_{\varepsilon}}

\newcommand{\lhs}{left hand side}

\DeclareMathOperator{\deriv}{d}

\DeclareMathOperator{\dist}{dist}

\let\TeXchi\chi
\def\chi{{\setbox0 \hbox{\mathsurround0pt
$\TeXchi$}\hbox{\raise\dp0 \copy0 }}}

\newcommand{\calX}{{\mathcal X}^\eb}

\newcommand{\calE}{{\mathcal E}}
\newcommand{\calEe}{{\mathcal E}_\epsi}

\newcommand{\calV}{{\mathcal V}^\eb}

\newcommand{\calY}{{\mathcal Y}}

\newcommand{\barO}{\overline{\Omega}}

\newcommand{\dit}{\deriv\!t}
\newcommand{\dis}{\deriv\!s}

\newcommand{\ddt}{\frac{\deriv\!{}}{\dit}}

\numberwithin{equation}{section}
\begin{document}
\title
{Trajectory and smooth attractors \\
for Cahn-Hilliard equations with inertial term}
\author{Maurizio Grasselli\\
{\sl Dipartimento di Matematica}\\
{\sl Politecnico di Milano}\\
{\sl Via Bonardi, 9}\\
{\sl I-20133 Milano, Italy}\\
{\rm E-mail:~~\tt maurizio.grasselli@polimi.it}
\and
Giulio Schimperna\\
{\sl Dipartimento di Matematica}\\
{\sl Universit\`a di Pavia}\\
{\sl Via Ferrata, 1}\\
{\sl I-27100 Pavia, Italy}\\
{\rm E-mail:~~\tt giusch04@unipv.it}
\and
Sergey Zelik\\
{\sl Department of Mathematics}\\
{\sl University of Surrey}\\
{\sl Guildford, GU2 7XH, United Kingdom}\\
{\rm E-mail:~~\tt S.Zelik@surrey.ac.uk}
}
\maketitle
\begin{abstract}

\noindent The paper is devoted to a modification of the classical
Cahn-Hilliard equation proposed by some physicists. This modification
is obtained by adding the second time derivative of the order
parameter multiplied by an inertial coefficient $\epsi>0$ which is
usually small in comparison to the other physical constants. The main
feature of this equation is the fact that even a globally bounded
nonlinearity is ``supercritical'' in the case of two and three space
dimensions. Thus the standard methods used for studying semilinear
hyperbolic equations are not very effective in the present case.
Nevertheless, we have recently proven the global existence and
dissipativity of strong solutions in the 2D case (with a cubic
controlled growth nonlinearity) and for the 3D case with small
$\epsi$ and arbitrary growth rate of the nonlinearity (see
\cite{GSZ,GSSZ}). The present contribution studies the long-time
behavior of rather weak (energy) solutions of that equation and it is
a natural complement of the results of our previous papers
\cite{GSZ} and \cite{GSSZ}. Namely, we prove here that the
attractors for energy and strong solutions coincide for both the cases
mentioned above. Thus, the energy solutions are asymptotically
smooth. In addition, we show that the non-smooth part of any
energy solution decays exponentially in time and deduce that the
(smooth) exponential attractor for the strong solutions constructed
previously is simultaneously the exponential attractor for the energy
solutions as well. It is worth noting that the uniqueness of energy
solutions in the 3D case is not known yet, so we have to use the
so-called trajectory approach which does not require the uniqueness.
Finally, we apply the obtained exponential regularization of the
energy solutions for verifying  the dissipativity of solutions of the 2D
modified Cahn-Hilliard equation in the intermediate phase space of
weak solutions (in between energy and strong solutions) without any
restriction on $\epsi$.
\end{abstract}

\noindent {\bf Key words:}\thspace trajectory attractors, smooth
global attractors, singularly perturbed Cahn-Hilliard equation.

\vspace{2mm}

\noindent
{\bf AMS (MOS) subject clas\-si\-fi\-ca\-tion:}
\thspace 35B40, 35B41, 82C26.

\vspace{2mm}



%
%


\section{Introduction}
\label{secintro}

In a series of contributions, P.~Galenko {\sl et al.} (see
\cite{GJ,GL1,GL2,GL3}) have proposed to modify the celebrated
Cahn-Hilliard equation (see \cite{CH}, cf. also the review \cite{N-C})
in order to account for nonequilibrium effects in spinodal
decomposition (cf. \cite{Ca}, see also \cite{Gr,MW}). The basic form
of this modification reads
\begin{equation}\label{CHin}
  \epsi u\dtt+u_t-\Delta(-\Delta u+f(u))=g,
\end{equation}
on $\Omega\times (0,+\infty)$, $\Omega$ being a bounded smooth
subset of $\RR^{\rm N}$, ${\rm N}\leq 3$. Here $\epsi\in (0,1]$,
$f$ is the derivative of a nonconvex potential (e.g., $f(r)=r(r^2-1)$)
and $g$ is a given (time-independent) function.

The longtime behavior of equation \eqref{CHin} already drew the
attention of mathematicians (see the pioneering \cite{D}, cf. also the
more recent \cite{BGM,GGMP1D,GGMP1DM,ZM1,ZM2}). However, all
these contributions were essentially devoted to the one-dimensional
case which can now be considered well known. There also have been
further works devoted to higher dimensions (see \cite{GGMP3D,Ka},
cf. also \cite{CM,Ve} for memory effects) but they are all
characterized by the presence of viscosity terms which imply the
instantaneous regularization of solutions. This is not the case of
\eqref{CHin}.

The mere existence of energy bounded solutions (see
Def.~\ref{reg-of-sols} for this terminology) was proven  in
\cite{Se07} for $N=3$. We recall that this work was mainly devoted
to the longtime behavior of such solutions based on the multi-valued
semigroup approach to the problems without uniqueness developed
by A.V.~Babin and M.I.~Vishik \cite{BV1}
(see also \cite{Ba}, \cite{ball1,ball2} and references therein). The existence result was then generalized to a nonisothermal
system with memory in \cite{LR}. Nevertheless, the existence of
energy bounded solutions is not a delicate issue and can be carried out
by means of a standard Galerkin procedure. In addition, both the
quoted results were proven supposing $f$ of cubic controlled growth,
but the existence also holds when  $f$ has a generic polynomial
growth (cf. Thm.~\ref{existener} below). On the contrary, uniqueness
of such solutions is much harder to prove. This was eventually done in
\cite{GSZ} for ${\rm N}=2$, assuming $f$ of cubic controlled
growth, along with a number of other results (e.g., existence of
smoother solutions, global attractors and exponential attractors).
Uniqueness of energy bounded solutions is still open in the case ${\rm
N}=3$ (and also when ${\rm N}=2$ for supercubic $f$). More recently,
an extension to a version with memory has been studied in \cite{CC}.
However, the existence of stronger solutions was shown in \cite{GSSZ} for
$\epsi$ small enough. This fact enabled the authors to construct a
dynamical system acting on a suitable phase space (depending on
$\epsi$) and to prove the existence of the global attractor as well as
of an exponential attractor.

Speaking of global attractors in the cases ${\rm N}=2$ and ${\rm
N}=3$, the results obtained in \cite{GSZ} and \cite{GSSZ} are not
fully satisfactory. Indeed, in the former case we proved the existence
of the global attractors both for energy bounded solutions and for
``quasi-strong" solutions (see Def.~\ref{reg-of-sols} again), but we
could not say whether they coincide. In the case ${\rm N}=3$ we
only established the existence of the global attractor for quasi-strong
solutions, while the unique available result on global attractors for
energy bounded solutions was in \cite{Se07}. Inspired by
\cite{Ze04}, here we intend to bridge this gap.
\par
First, in Sections \ref{secmath} and \ref{sectraj}, we construct the proper attractor
 for the energy solutions of \eqref{CHin}. Since we
do not have the uniqueness for the energy solutions (in the 3D case
as well as in the 2D case with the super-cubic growth rate), we use
the so-called trajectory dynamical system approach developed by
V.V.~Chepyzhov and M.I.~Vishik (see \cite{CV})  and construct the
so-called trajectory attractor associated with energy solutions of
problem \eqref{CHin}. However, the class of all energy solutions
which
 satisfy the weakened form of energy inequality used in
 \cite{CV} in their construction of trajectory attractors for
 damped hyperbolic equations is too large for our purposes and
 we restrict ourselves to consider only the energy solutions
 which can be obtained by Galerkin approximations
 (analogously to \cite{Ze04}, see also \cite{Se07}). Note that, although the trajectory attractor constructed here is very close (and even formally equivalent) to the generalized attractor obtained in \cite{Se07} via the multi-valued approach, it is much more convenient for our further investigation.
 \par
 Then, in Section \ref{secsmooth}, we establish that each complete bounded solution belonging to the trajectory attractor
is a strong solution to \eqref{CHin} at least for all times smaller than
a time sufficiently close to $-\infty$. This kind of backward
smoothness was firstly obtained in \cite{Ze04} for damped wave
equations with supercritical nonlinearities. Here, the proof is however
based on  partly different and more simple arguments (see Thm.~\ref{Zelik}).
\par
 The backward smoothness is the basic ingredient which allows us to show (in Section \ref{smoothattr}):

\smallskip\noindent
(i) if ${\rm N}=2$ (and $f$ has cubic controlled growth) the global
attractor for energy bounded solutions
    coincides with the one for quasi-strong solutions (and, in particular,  it is smooth);

\smallskip\noindent
(ii) if ${\rm N}=3$ the trajectory attractor consists of
    complete bounded strong solutions if $\epsi$ is small enough.

\smallskip\noindent
Thus, in both cases, we have the asymptotic regularization of energy solutions.
 \par
In Section \ref{s5.exp}, in both cases mentioned above, we establish
that every energy solution regularizes {\it exponentially} as $t\to
+\infty$, i.e., such solutions can be split into the sum of two
functions, one of which is smooth and bounded and the other tends
to zero exponentially as time tends to infinity. This result seems new
even for the well-known damped semi-linear wave equation with
supercritical nonlinearity (a similar property has been shown in
\cite{Ze04} {\it only} under the additional assumptions that all
equilibria are hyperbolic). In addition, we present an alternative
approach to demonstrate the exponential asymptotic regularization
when $\eb>0$ is small enough. This method does not use the
backward regularization or exploit the global Lyapunov functional and can
be therefore applied, e.g., to non-autonomous equations or to the
case of unbounded domains.
\par
Finally, in Section \ref{s.exp}, taking advantage of the exponential regularization and the transitivity of exponential attraction, we prove that, again in both the above cases, the energy solutions approach {\it exponentially} fast the exponential attractor for strong solutions which has been constructed in \cite{GSSZ} and \cite{GSZ} . Thus, this strong exponential attractor is
the exponential attractor for the energy solutions as well. In addition, we use the obtained exponential regularization to solve one problem for the 2D case which remained open in \cite{GSZ}, namely, the dissipativity in the intermediate phase space of weak solutions (between the energy and strong solutions) with no restrictions on $\eb$.
\par
 To conclude the introduction, we note that, although
 we endow equation \eqref{CHin}  with the
boundary and initial conditions
\begin{align}
  &\disp u(t) = \Delta u(t)=0, \quad \mbox{ on } \partial\Omega, \; t>0 ,\label{bc} \\
  &\disp u(0) = u_0, \quad u_t(0)= u_1, \quad  \mbox{ in } \Omega,
  \label{ic}
\end{align}
other boundary conditions, like no-flux or periodic, could also be
handled (see \cite{BGM,D,GGMP1D}) with only technical modifications.


\section{Functional setup and existence of solutions}
\label{secmath}

Let us set $H:=L^{2}(\Omega)$ and denote by
$(\cdot,\cdot)$ the scalar product both in $H$ and
in $H\times H$, and by $\|\cdot\|$ the related norm.
The symbol $\|\cdot\|_{X}$ will indicate the norm
in the generic (real) Banach space $X$.
Next, we set $V:=H^1_0(\Omega)$, so that
$V'=H^{-1}(\Omega)$ is the topological
dual of $V$. The duality between $V'$ and $V$ will
be noted by $\langle\cdot,\cdot\rangle$.
The space $V$ is endowed with
the scalar product
\begin{equation} \label{proscaV}
  (\!(v,z)\!):=\io\nabla v\cdot\nabla z, \quad\perogni v,z\in V,
\end{equation}
and the corresponding induced norm. We shall denote by $A$ the
Riesz operator on $V$ associated with the norm above, namely,
\begin{equation} \label{AVVp}
  A:V\to V', \qquad
   \duav{Av,z}=(\!(v,z)\!)=\io\nabla v\cdot\nabla z,
    \quad \perogni v,z\in V.
\end{equation}
%
%
%
Abusing notation slightly, we shall also indicate
by the same letter $A$ the restriction of
the operator defined in \eqref{AVVp} to the set
$D(A)=H^2(\Omega)\cap V$, i.e., the unbounded operator defined as
\begin{eqnarray}\label{defA}
  A=-\Delta \;\;\;\;\mbox{ with domain }\,
   D(A)=H^{2}(\Omega)\cap V\subset L^{2}(\Omega).
\end{eqnarray}
Starting from $A$ one can define the family of Hilbert spaces
\begin{equation*}
  H^{2s}=D(A^{s}), \quad \;s\in \mathbb{R},\\
\end{equation*}
with scalar product $(A^{s}\cdot,A^{s}\cdot)$.
It is well known that $H^{s_1}\subset H^{s_2}$
with dense and compact immersion when $s_1>s_2$.
Then, we introduce the scale of Hilbert spaces
\begin{equation}
  \calV_s:=
    D(A^{\frac{s+1}2})\times \sqrt{\epsi}D(A^{\frac{s-1}2}),
\end{equation}
so that we have, in particular, $\calV_0=V\times V'$ and,  analogously,
$\calV_1=(H^2(\Omega)\cap V)\times H$.
The spaces $\calV_s$ are
naturally endowed with the graph norm
\begin{equation}\label{definorVs}
  \|(u,v)\|_{\calV_s}^2 := \|A^{\frac{s+1}2} u\|_H^2
   +\epsi\|A^{\frac{s-1}2}v\|_H^2.
\end{equation}

Regarding the nonlinear function $f$, we assume that
$f\in C^3(\RR;\RR)$ with $f(0)=0$ satisfies, for some $p\in[0,\infty)$,
\begin{align}
  & \disp \liminf_{\vert r\vert\nearrow +\infty}\frac{f(r)}{r}> -\lambda_1;
 \label{f1}\\
  & \disp \esiste \lambda \in [0,+\infty) \mbox{ and } \delta\in[0,\infty):~~
   f'(r)\ge-\lambda+\delta|r|^{p+2},~~\perogni r\in\RR;\label{f2}\\
 %
 %
 \label{f3}
  & \disp\esiste M\ge 0:~~|f'''(r)|\le M(1+|r|^p), ~~\perogni r\in\RR.
\end{align}
Here $\lambda_1>0$ is the first eigenvalue of $A$.
Note that $f$ can be a polynomial of arbitrarily large odd degree
with positive leading coefficient.
If we indicate as $F $ the potential of $f$ (i.e.,
a suitable primitive of $f$), we can always suppose that
\begin{equation}\label{prim}
  F(r) \ge -\frac{\kappa}{2}r^2,
\end{equation}
for some $\kappa<\lambda_1$. By \eqref{f2},
we also have that $F$ is $\lambda$-convex. When $p>2$
we will need to suppose $\delta>0$ in \eqref{f2}.

Finally, we let
\begin{equation}
  g \in H.\label{regog}
\end{equation}
%
%
%
%
System \eqref{CHin}-\eqref{ic} can then be reformulated as

\medskip
\noindent
\textbf{Problem $P_\epsi$.} \label{problema debole}
Find  a pair $(u,u_t)$ 
satisfying
\begin{align}
  &\disp \epsi u_{tt} + u_{t} +A(Au +f(u)) = g, 
      \label{CH1-regular} \\
  &\disp u|_{t=0}=u_{0}, \quad u_{t}|_{t=0}=u_{1}. \label{iniz-regular}
\end{align}
Here, it is intended that the first two relations
hold at least for almost any time $t$ in
the life span of the solution. Actually,
we will consider in the sequel both local and global
in time solutions. In the sequel we will
frequently write $U$ for the couple $(u,u_t)$
and $U_0$ for $(u_0,u_1)$ (the same convention
will be kept for other letters, e.g., we will write
$V=(v,v_t)$). Moreover, for the sake of brevity,
solutions will be sometimes noted simply as $u$,
or as $U$, rather than as $(u,u_t)$.

Speaking of regularity, we can now introduce
the {\sl energy}\/ associated with \eqref{CHin} as
\begin{equation}\label{defiEe}
  \calEe:\calV_0\to \RR, \qquad
   \calEe(u,v):=\frac12\|(u,v)\|_0^2+\io F(u)
     -\duav{g,A^{-1}u}.
\end{equation}
%
%
%
Assumptions \eqref{f3} and \eqref{regog} suffice to guarantee that
$\calEe$ is finite for all $(u,v)\in \calV_1$. However, if $p>2$
then $\delta>0$ in \eqref{f2} is needed, if $(u,v)\in \calV_0$.
For this reason we introduce the function space
\begin{equation}\label{defiXe0}
  \calX_0:=\big\{(u,v)\in\calV_0:
   u\in L^{p+4}(\Omega)\big\},
\end{equation}
which is endowed with the graph metrics.
For instance, with some abuse of language,
we will write
\begin{equation}\label{defimetricXe0}
  \|(u,v)\|_{\calX_0}^2:=
   \|(u,v)\|_{\calV_0}^2+ \|u\|^{p+4}_{L^{p+4}(\Omega)}.
\end{equation}
It is then clear from \eqref{f3} that
$\calEe$  is locally finite in $\calX_0$.
Of course, thanks to \eqref{f1}, $\calEe$
is in any case bounded from below on the whole
$\calV_0$.
The above discussion leads to the following
definition (see \cite{GSSZ,GSZ}).
\bede\label{reg-of-sols}
 We say that a solution to~$P_\epsi$ defined on some time interval $(0,T)$
 is an\/ {\rm energy bounded solution}, or, more concisely,
 {\rm energy solution}, if
 $(u,u_t)\in L^{\infty}(0,T;\calV_{0})$
 and $\calE_\epsi(u,u_t)\in L^\infty(0,T)$.
 If $(u,u_t)\in L^{\infty}(0,T;\calV_1)$, we say instead that $u$ is a\/ {\rm weak
 solution}.
 If $(u,u_t)\in L^{\infty}(0,T;\calV_2)$, $u$ is named \/ {\rm quasi-strong
 solution}, while $u$ is a \/ {\rm strong
 solution} if $(u,u_t)\in L^{\infty}(0,T;\calV_3)$.
\edde
Thus, for energy solutions, \eqref{CH1-regular} has to be
interpreted as an equation in $D(A^{-2})$ in the case when $p>2$
(and hence $\delta>0$). Indeed, in this case $f(u(t))\in
L^{\frac{p+4}{p+3}}(\Omega)\subset D(A^{-1})$ for almost any
$t\in (0,T)$. If $p\le 2$, then of course we can say more: it is now
$f(u(t))\in L^{6/5}(\Omega)\subset V'$ and \eqref{CH1-regular}
holds in $D(A^{-3/2})$. In any case, we can say that  any energy
solution $U=(u,u_t)$ lies in $L^{\infty}(0,T;\calX_0)$. Passing to
weak solutions, then \eqref{CH1-regular} holds in $D(A^{-1})$
since it is now $f(u(t))\in H$, thanks to the embedding
$H^2(\Omega)\hookrightarrow C(\barO)$. For the same reason,
\eqref{CH1-regular} can be interpreted as a $V'$-equation for
quasi-strong solutions, and, of course, for strong solutions it holds
almost everywhere in $\Omega\times(0,T)$. Observe that, despite
of the name, energy bounded solutions are {\sl weaker}\/ than
weak solutions.

We also notice that a comparison in \eqref{CH1-regular} gives that
$u_{tt}\in L^{\infty}(0,T,H)$ for strong solutions,
$u_{tt}\in L^{\infty}(0,T,V')$ for quasi-strong solutions,
$u_{tt}\in L^{\infty}(0,T,D(A^{-1}))$ for weak solutions,
and $u_{tt}\in L^{\infty}(0,T,D(A^{-2}))$ for energy
solutions. This immediately leads to
$U=(u,u_t)\in C^{0}_{w}([0,T];\calV_i)$ with,
respectively, $i=3,2,1,0$, where
$C^{0}_{w}([0,T];X)$ is defined as ($X$ being a real Banach space)
$$
  C^{0}_{w}([0,T];X):=\left\{v\in L^{\infty}(0,T;X): \;\langle \phi,v(\cdot)\rangle
   \;\in C^0([0,T]), \;\;\forall \phi\in X'\right\}.
$$
Therefore solutions can be evaluated pointwise in
time and initial conditions \eqref{iniz-regular}
in $\calV_i$, $i=3,2,1,0$, have a well-defined meaning in all the cases.

We conclude the section by stating the existence theorem
\bete
\label{existener}
 Let the assumptions \eqref{f1}-\eqref{f3}
 and \eqref{regog} hold,
 and let
 \begin{equation}\label{cond-init-w}
   (u_0,u_1)\in \mathcal{X}_0.
 \end{equation}
 Then, if either $p\le 2$ in \eqref{f2}-\eqref{f3}
 or $\delta>0$ in \eqref{f2},
 there exists at least one\/ {\rm global in time}
 energy solution to
 Problem~$P_\epsi$, which additionally satisfies
 the following dissipation inequality
 \begin{equation}\label{weakenergy}
   \disp\|U(t)\|^2_{\calX_0} + \int_{t}^{+\infty} \|u_t(s)\|^2_{V'}\,\dis
   \disp \le C\|U_0\|^2_{\calX_0}e^{-\kappa t} + C(1+\|g\|^2),
   \quad\forall t\ge \tau\ge 0,
 \end{equation}
 where the positive constants $\kappa$ and $C$ are independent of $\epsi$.
\ente
From now on we let the assumptions of Theorem \ref{existener}
hold, unless otherwise specified. Moreover, we will not stress the
dependence on $\eps$ till the final section.


Although the proof of Theorem~\ref{existener} is standard, we report
here below some highlights for the reader's convenience.
\begin{proof}
 The proof is essentially based on the Faedo-Galerkin scheme described
 in the next section and on a couple of estimates which, for simplicity,
 are performed here by working directly (albeit formally) on the
 original problem rather than on its approximation. Firstly,
 testing \eqref{CH1-regular} by $A^{-1}u_t$, we easily derive
 the {\sl energy equality}
 \begin{equation}\label{en-eq}
   \ddt \calE\ee(U)+\|u_t\|_{V'}^2=0
 \end{equation}
 (in fact, the above is a true equality just at the regularized
 (Galerkin) level, but will turn into an inequality
 when taking the limit). Second,
 multiplying \eqref{CH1-regular} by $A^{-1}u$, we get
 \begin{equation}\label{sec-eq}
   \ddt \Big[ \epsi\duav{u_t,A^{-1}u} + \frac12\|u\|_{V'}^2 \Big]
    - \epsi  \| u_t \|_{V'}^2 + \| u \|_V^2 + (f(u),u) - (g,A^{-1}u) = 0.
 \end{equation}
 We then notice that a combination of assumptions
 \eqref{f1}-\eqref{f2} gives
 \begin{equation}\label{giulio1}
   \| u \|_V^2 + (f(u),u)
    \ge \kappa_1\| u \|_V^2 + \kappa_2\io F(u) - c
    \ge \kappa_3\| u \|_V^2 - c,
 \end{equation}
 for suitable positive constants $\kappa_i$, $i=1,2,3$, only
 depending on $\lambda_1$ and $\lambda$.
 In particular, thanks to \eqref{f1},
 $\kappa_2$ can be chosen so small that
 the latter inequality holds even in case $\delta=0$ in
 \eqref{f2}.

 Using \eqref{giulio1},
 it is a standard matter to verify that, multiplying
 \eqref{sec-eq} by a (suitably small) constant $\alpha>0$,
 and adding the result to \eqref{en-eq}, gives,
 for some $\kappa>0$,
 \begin{equation}\label{diss-eq}
   \ddt \calY\ee + \kappa \calY\ee
     \le c,
 \end{equation}
 where we noted as $\calY\ee$ the functional obtained
 by taking $\calE\ee$ and summing to it $\alpha$ times
 the terms in square brackets in \eqref{sec-eq}
 and adding also a further quantity
 of the form $C_*(1+\|g\|^2)$, where
 $C_*>0$ is a suitable constant. Actually, by
 \eqref{defimetricXe0}, it is easy to check
 that, if $C_*$ is taken large enough, then
 for some positive $\kappa_i$, $i=4,\dots,7$,
 and $c>0$ there holds
 \begin{equation}\label{prop-Y}
   \kappa_4 \| U \|_{\calV_0}^2 + \delta\|u\|^{p+4}_{L^{p+4}(\Omega)}
   \le \kappa_5\calE\ee
    \le \calY\ee
    \le \kappa_6\calE\ee + c (1+\|g\|^2_H)
    \le \kappa_7\| U \|_{\calX_0}^2 + c (1+\|g\|_H^2).
 \end{equation}
 Thus, integrating \eqref{diss-eq} over $(0,t)$ and
 using \eqref{prop-Y}, we obtain an inequality analogue
 to \eqref{weakenergy}, but without the
 integral term on the \lhs. To get a control of
 it, it is however sufficient to go back
 to \eqref{en-eq}, integrate it over $(t,+\infty)$,
 and refer once more to \eqref{prop-Y}.

 Once one has obtained \eqref{weakenergy} at the
 approximated level, it is then a standard procedure
 to pass to the limit. For instance, assumptions
 \eqref{f2}-\eqref{f3} and the Aubin-Lions
 compactness lemma entail weak star convergence
 of $f(u)$ to the right limit, holding
 in the space $L^\infty(0,T;L^{6/5}(\Omega))$ if
 $p\le 2$, and in
 $L^\infty(0,T;L^{\frac{p+4}{p+3}}(\Omega))$
 if $p>2$ (and hence $\delta>0$). Then,
 also \eqref{weakenergy} passes to the  limit inferior
 by standard semicontinuity arguments. The proof is complete.
\end{proof}


\begin{osse}
\label{weakf}
It is not difficult to realize that Theorem~\ref{existener} still holds when
$f\in C^0(\RR;\RR)$ with $f(0)=0$ satisfies, in place of \eqref{f1}-\eqref{f3} and
for some $p\in(0,\infty)$,
\begin{align*}
& \disp\esiste M_0\ge 0:~~|f(r)|\le M_0(1+|r|^{p+3}), ~~\perogni r\in\RR\\
  & \disp \esiste \lambda_0 \in [0,+\infty) \mbox{ and } \delta_0\in (0,\infty):~~
   f(r)r\ge-\lambda_0+\delta_0|r|^{p+4},~~\perogni r\in\RR.
\end{align*}
If $p\le2$, the last condition can be replaced by \eqref{f1}.
These assumptions on $f$ also suffice to establish the existence of the trajectory
attractor (see Theorem~\ref{Zelik0} and related corollaries in the next section).
\end{osse}


\section{The trajectory dynamical system}
\label{sectraj}

The existence of a global energy solution (see Theorem \ref{existener}) can
be proven by means of a Faedo-Galerkin procedure similar to the one used in
\cite{Ze04} for the damped semilinear wave equation with a supercritical
nonlinearity. In particular, if we indicate by $P_n$ the orthoprojector
constructed with the first $n$ eigenfunctions of $A$ and we set
\begin{equation}
U^n_0=(u_0^n,u_1^n)=(P_n u_0,P_n u_1) \in \mathcal{V}_{0n}:=
(P_n H)^2,
\end{equation}
then it is not difficult to prove that the corresponding approximating solution
$U^n(t) = (u^n(t),u^n_t(t))$ to
\begin{align}
  &\disp \epsi u^n_{tt} + u^n_{t} +A(Au^n +P_n f(u^n)) = g^n:= P_n g, 
      \label{CH1-n} \\
  &\disp u^n|_{t=0}=u^n_{0}, \quad u^n_{t}|_{t=0}=u^n_{1}, \label{iniz-n}
\end{align}
also satisfies the analogue of \eqref{weakenergy}, namely,
\begin{equation}\label{weakenergyn}
   \disp\|U^n(t)\|^2_{\calX_0} + \int_{t}^{\infty} \|u^n_t(s)\|^2_{V'}\,\dis
   \disp \le C\|U^n(\tau)\|^2_{\calX_0}e^{-\kappa(t-\tau)} + C(1+\|g\|_H^2),
   \quad\forall t\ge\tau\ge 0.
\end{equation}
From now on, we restrict ourselves to consider only those of energy
bounded solutions which can be obtained as a weak limit of the
corresponding Galerkin approximations. \bede\label{Def.ensol} An
energy bounded solution $U(t):=(u(t),u_t(t))$ of problem
\eqref{CHin} is an {\it energy} solution if it can be obtained as a
weak limit of a subsequence of solutions $U_n(t)$ to the Galerkin
approximation equations \eqref{CH1-n} which satisfy
\eqref{iniz-n}. \edde
Since the uniqueness of the energy solutions is
not known so far, we use the so-called trajectory approach
developed in \cite{CV} in order to describe the long-time behavior of
such solutions. To this end, we first need to define the trajectory
phase space associated with energy solutions of problem
\eqref{CHin} and the trajectory dynamical system on it.

\bede\label{Def.tra} We define the trajectory phase space $K^+_\eb\in L^\infty(\R_+,\calX_0)$ as a set of
all energy solutions $U$ of problem \eqref{CHin} associated with all possible initial data $U(0)\in\calX_0$. Namely,
we set
\begin{align}
K^+_\eb := \Big\{&U\in L^{\infty}({\mathbb R}_+;\calX_0)\; : \; \exists \{U^{n_k}(t)\}
\textrm{ solving } \eqref{CH1-n}-\eqref{iniz-n} \textrm{ such that }\\
\nonumber
&U(0) = [\calX_0]^w - \lim_{k\to \infty} U^{n_k}(0) \textrm{ and }
U = \Theta^+ - \lim_{k\to \infty}U^{n_k}\Big\}
\end{align}
and we endow $K^+_\eb$ with the topology of $\Theta^+:=[L^\infty_{loc}([0,\infty),\calX_0)]^{w^*}$ (the weak star topology  of $L^\infty_{loc}([0,\infty),\calX_0)$).
Then, the time translation semigroup
\begin{equation}
\mathbb{T}_\ell: K^+_\eb \to K_\eb^+, \qquad (\mathbb{T}_\ell u)(t) = u(t+\ell),
\end{equation}
is well defined for $\ell\ge 0$. The semigroup  $\mathbb{T}_\ell$ acting on $K^+_\eb$ endowed by the above defined topology is called the trajectory dynamical system associated with equation \eqref{CHin}.
\edde
We recall that a sequence $\{V^n\}\subset
L^{\infty}({\mathbb R}_+;\calX_0)$ converges to $V$ in
$\Theta^+$ if, for every $T \ge 0$, $V^n \to V$ weakly star in
$L^\infty((T,T+1);\calX_0)$ as $n$ goes to $\infty$. Similarly, we
can endow $L^{\infty}({\mathbb R};\calX_0)$ with the weak star
local topology $\Theta:=[L^\infty_{loc}(\R,\calX_0)]^{w^*}$. The obtained topological
spaces are Hausdorff and Fr\'{e}chet-Urysohn with a countable base of
open sets (see \cite[Chap.~XII]{CV}).
\par
In order to be able to speak about the attractor of the trajectory dynamical system $(\mathbb T_\ell,K^+_\eb)$ (i.e., the trajectory attractor of equation \eqref{CHin}), we also need to define the class of {\it bounded} sets in a proper way. We note that,
in contrast to \cite{CV}, only the energy bounded solutions which can be obtained through the Galerkin limit are included in
$K^+_\eb$ (and that difference is crucial for what follows). By this reason, we need, in addition, to introduce (following \cite{Ze04}) the so-called $M$-functional on the space $K^+_\eb$:
\begin{equation}
\label{Mfunctional}
M_u^\eb(t) := \inf \Big\{\liminf_{k\to \infty} \Vert U^{n_k}(t)\Vert_{\calX_0}
\; : \;
U  = \Theta^+ - \lim_{k\to \infty}U^{n_k}, \;
U(0) = [\calX_0]^w  - \lim_{k\to \infty}U^{n_k}(0)
\Big\},
\end{equation}
where the infimum is taken over all the sequences $\{U^{n_k}(t)\}_{k\in {\mathbb N}}$ of
Faedo-Galerkin approximations which
$\Theta^+$-converge to the given solution $U$.

Recalling now \cite[Cor. 1.1]{Ze04}, we can easily prove the following properties
of the $M-$energy functional, i.e., for any $U\in K^+_\eb$ we have
\begin{align}
\label{M1}
&M^\eb_u(t) < \infty, \quad \Vert U(t) \Vert_{\calX_0} \leq M^\eb_u(t), \quad
  M^\eb_{\mathbb{T}_\ell u}(t)\leq M^\eb_{u}(t+\ell),\\
\label{M2}
&M^\eb_u(t)^2 +  \int_{t}^{\infty} \|u_t(s)\|^2_{V'}ds\,
\le CM^\eb_u(\tau)^2e^{-\kappa(t-\tau)} + C_0(1+\|g\|_H^2),
   \quad\forall t\ge\tau\ge 0.
\end{align}
We can now say that a set $B\subset K^+_\eb $ is $M-$bounded if
\begin{equation}
\sup_{U\in B} M^\eb_u(0) < \infty
\end{equation}
and recall the definition of the trajectory attractor associated with  \eqref{CHin}.
\begin{defi}\label{Def.tr} A set $\mathcal A_\eb^{tr}\subset K^+_\eb$ is a trajectory attractor associated with energy solutions of equation \eqref{CHin} (i.e., the global attractor of the trajectory dynamical system $(\mathbb T_\ell,K^+_\eb)$) if:
\par
I) the set $\mathcal A^{tr}_\eb$ is compact in $K^+_\eb$ (endowed by the $\Theta^+$ topology);
\par
II) it is strictly invariant: $\mathbb T_\ell\mathcal A^{tr}_\eb=\mathcal A^{tr}_\eb$, $\ell\ge0$;
\par
III) for every $M$-bounded set $B\subset K^+_\eb$ and every neighborhood $\mathcal O(\mathcal A^{tr}_\eb)$ of $\mathcal A^{tr}_\eb$
(again in the topology of $\Theta^+$), there exists $T=T(B,\mathcal O)$ such that (attraction property)
$$
\mathbb T_\ell B\subset\mathcal O(\mathcal A^{tr}_\eb),\quad \forall \ell\ge T.
$$

\end{defi}
We can now state the existence of the trajectory attractor which can be proven
arguing as in \cite[Thm.~1.1]{Ze04}.

\bete \label{Zelik0}  Let \eqref{f1}-\eqref{f3} and \eqref{regog}
hold. Then the semigroup $\mathbb{T}_\ell$ acting on $K^+_\eb$
possesses the trajectory attractor ${\mathcal A}^{tr}_\eb$
characterized as follows
\begin{equation}
\label{trajattr}
{\mathcal A}^{tr}_\eb = \Pi_+ {\mathcal K}_\eb,
\end{equation}
where ${\mathcal K_\eb}\in L^\infty({\mathbb R};\calX_0)$ is the set of
all the complete $\calX_0$-bounded solutions to $P_\epsi$ which can be obtained as a Faedo-Galerkin
limit, while $\Pi_+$ is the projection onto $L^\infty({\mathbb R}_+;\calX_0)$.
More precisely, $ U \in {\mathcal K_\eb}$ if and only if there exist $\{t_k\}$
such that $t_k \searrow -\infty$ and $\{U^{n_k}(t)\}$ such that, for $t\geq t_k$,
\begin{align}
  &\disp \epsi u^{n_k}_{tt} + u^{n_k}_{t} +A(Au^{n_k} +P_{n_k} f(u^{n_k}))
  = g^{n_k}, 
      \label{CH1-nk} \\
  &\disp u^{n_k}|_{t=t_k}=u^k_{0}, \quad u^{n_k}_{t}|_{t=t_k}=u^k_{1}, \label{iniz-nk}
\end{align}
with
\begin{equation}
\label{bdd-conv-nk}
U^{n_k}(t_k)\in (P_{n_k}H)^2,\qquad
\Vert U^{n_k}(t_k)\Vert_{\calX_0} \leq C,\qquad
U = \Theta - \lim_{k\nearrow \infty} U^{n_k},
\end{equation}
where $C>0$ is independent of $k$.
\ente
The proof of this theorem repeats word by word the proof of
\cite[Thm.~1.1]{Ze04} and, for this reason, is omitted.
Still arguing as in \cite{Ze04}, we can deduce the following corollaries.
In the statements, if $X,Y$ are Banach spaces, $C_b(X;Y)$ will denote
the Banach space of all continuous {\sl and bounded}\/ functions
from $X$ to $Y$, endowed with the supremum norm.
\beco
\label{attraction}
Let $B\subset K^+$
an M-bounded set. Then, for every $T\in {\mathbb R}_+$ and every $\beta \in (0,1]$,
the following convergence holds
\begin{equation}\label{conv}
\lim_{\ell \to \infty}
\dist_{{\mathcal L}_\beta(\ell,T+\ell)}(B\vert_{[\ell,T+\ell]}, {\mathcal A}^{tr}\vert_{[\ell,T+\ell]})=0,
\end{equation}
where
\begin{equation}
{\mathcal L}_\beta(\ell,T+\ell)= C([\ell,T+\ell],[D(A^{(1-\beta)/2})\cap L^{p+4-\beta}(\Omega)]
\times \sqrt{\epsi}D(A^{-(1+\beta)/2})).
\end{equation}
Here we recall the definition of the Hausdorff semidistance, namely
\begin{equation}
\dist_{\mathcal{L}}(\mathcal{B}_1, \mathcal{B}_2)
:= \sup_{u\in \mathcal{B}_1} \inf_{v\in \mathcal{B}_2} d_{\mathcal{L}}(u,v),
\end{equation}
where $\mathcal{L}$ is some given metric space with distance $d_{\mathcal{L}}$ and
$\mathcal{B}_j\subset \mathcal{L},\,j=1,2$.
\enco

\beco
\label{dissipation}
Let $U\in {\mathcal K}_\eb$. Then, we have
\begin{equation}
\label{dissint}
\int_{-\infty}^{\infty} \Vert u_t(s) \Vert^2_{V'} ds\leq c(1 + \Vert g \Vert_H^2),
\quad u_{tt} \in C_b({\mathbb R},D(A^{-2})).
\end{equation}
Thus, for every $\beta >0$, there hold
\begin{equation}
\label{zelik1}
u_t \in C_b({\mathbb R},D(A^{-(1+\beta)/2})),  \quad
\lim_{t\to \pm \infty}\Vert A^{-(1+\beta)/2} u_t(t) \Vert_H^2 =0,
\end{equation}
and we also have the convergence to the set of equilibria ${\mathcal R}$
\begin{equation}
\label{zelik2}
\dist_{{\left(D(A^{(1-\beta)/2})\cap L^{p+4-\beta}(\Omega)\right)
\times \sqrt{\epsi} D(A^{-(1+\beta)/2})}}(U(t),{\mathcal R})\to 0,
\end{equation}
as $t$ goes to $\infty$, for any $\beta\in (0,1]$.
\enco

\begin{osse}
\label{choices} The choice of defining a trajectory dynamical
system by selecting those solutions which are limits of Galerkin
approximations excludes other possible solutions $(u,u_t)\in
L^\infty(\R_+,\calX_0)$ which satisfy the equation \eqref{CHin}
in the sense of distributions but cannot be obtained in that way.
Actually, nothing is known about such ``pathological'' solutions and
theoretically they may exist and even may not be dissipative (that is,
they may not satisfy the energy inequality). We remind that we cannot exclude that
this might happen even in dimension two when $f$ is
allowed to have a supercubical growth. Thus, in order to be able to
deal with attractors (no matter using the trajectory or multi-valued
semigroup approaches), one should restrict the admissible set of
solutions. To do that, there are at least two alternative ways.  The
first one is to consider only the solutions which satisfy some
weakened form of energy inequality, like the 3D Navier-Stokes equations
(see, e.g., \cite{CV} and references therein). The second one (used in
\cite{Ze04} and in this paper) is to consider only the solutions
which can be obtained by the Galerkin approximations (see
\cite[Rem.~6.2]{MZ} for more details). We only mention here that
both of them present some drawbacks. In particular, the class of
energy solutions may depend on some artificial constants
in the energy inequality in the former approach,
while it may depend on the choice of a Galerkin basis in the latter
one. Nonetheless, the second approach has an advantage which is
crucial in the present case, namely, it allows to justify the further
energy-like inequalities for the energy solutions. In fact, it is
completely unclear how to verify most of the results of this paper
using the first method. One more drawback is that both the
mentioned approaches destroy the concatenation property of energy
solutions and, up to the moment, no reasonable way to preserve this
property and exclude the ``pathological'' non-dissipative solutions is
known. Thus, the concatenation property seems to be an extremely
restrictive assumption which can be verified only in relatively simple
cases (usually, when the non-uniqueness is simply due to the
presence of non-Lipschitz nonlinearities). For this reason, in the
multi-valued approach, one usually needs (following \cite{BV1}, see
also \cite{RSS} and its references) to replace the semigroup {\it identity}
by the semigroup {\it inclusion} and the constructed attractor will
be also only {\it semi}-invariant with respect to the semigroup
itself. This is exactly the approach applied in \cite{Se07} to analyze
equation \eqref{CHin}. However, even though the trajectory
approach is formally equivalent, the trajectory attractor remains
strictly invariant even without the concatenation property (see
Definition~\ref{Def.tr}) and that makes it more convenient for
concrete applications.
\end{osse}


\section{Backward smoothness of complete trajectories}
\label{secsmooth}

In this section, we prove that any $U=(u,u_t)\in {\mathcal K}_\eb$ (see
Theorem \ref{Zelik0}) is backward smooth, i.e., $U(t)\in \calV_3$ if
$t$ is small enough. This kind of result is similar to
\cite[Thm.~2.1]{Ze04}; however, here we follow an alternative
strategy which makes use of stationary solutions. Actually, in
\cite{Ze04}, the globally defined trajectory was compared with the
corresponding trajectory of the limit evolution equation obtained by
setting $\varepsilon=0$. Here, instead, we first find a solution $v$ to
an auxiliary equation and we prove that $v$ is backward smooth.
Then, we show that actually $v(t)\equiv u(t)$ for all $t\leq T$, when
$T$ is small enough. A crucial point of the argument is the
construction of a smooth solution which will then be viewed as the
nonvanishing part of $v$. This is the content of the following

\bete\label{Zelik}
Let $U=(u,u_t)\in \mathcal{K}_\eb$. Then, for every
$\sigma>0$, there exist $T=T(\sigma,u)<0$ and a function $\tilde
u=\tilde u_\sigma \in C^\infty({\mathbb R}_-; H^4)$ such that:

\item{1)} for every $\beta >0$, there exists $c>0$ independent of $\sigma$ such that,
for every $t\le T$,
\begin{equation}
\label{Z1}
\|u(t)-\widetilde u(t)\|_{H^{1-\beta}}+\|u_t(t)-\widetilde u_t(t)\|_{H^{-\beta}}
\le c\sigma;
\end{equation}
\item{2)} there exists $C>0$, which is independent of $\sigma$ and $t$, such that
\begin{equation}
\label{Z4}
\|\tilde u\|_{C^2({\mathbb R}_-; H^4)}\le C(1+\Vert g \Vert_H);
\end{equation}
\item{3)} for each $m\in {\mathbb N}$, $m\geq 1$, there exist $C_m>0$
and $\kappa>0$, which are independent of $\sigma$ and $t$, such that
\begin{equation}
\label{Z5}
\|\partial^m_t \tilde u(t)\|_{H^4} \le C_m\sigma^{\kappa},
\qquad \forall\, t\le T.
\end{equation}
\item{4)} $\tilde u$ solves
\begin{equation}
\label{Z2}
\varepsilon\tilde u_{tt}+ \tilde u_t +A^2\tilde u+Af(\tilde u)=g +\phi(t),
\end{equation}
where $\phi$ is a suitable function such that
\begin{equation}
\label{Z3}
\|\phi(t)\|_H + \|\phi_t(t)\|_H \le C\sigma^{\kappa},\qquad \forall\, t\le T.
\end{equation}
\ente

\begin{proof}
Using \eqref{zelik1} and \eqref{zelik2}, we see that, for every
$\sigma>0$ and every $S\in {\mathbb N}$, there exists $T=T(\sigma,u,S)<0$ such
that, for every $s\le T$, there is an equilibrium $u_s$ which satisfies, for
some $\beta\in (0,1]$,
\begin{equation}
\label{zelik3}
\sup_{t\in[s,s+S]}\|u(t)-u_s\|_{H^{1-\beta}}\le \sigma.
\end{equation}
In fact, fixing say $S=2$ is sufficient for the proof.

Let us check now that
\begin{equation}
\label{zelik4}
\|u_{s+S/2}-u_s\|_{H^4} \le C\sigma^\kappa
\end{equation}
for some positive $C$ and $\kappa$ independent of $\sigma$ and $s$.
Indeed, since $u_s$ and $u_{s+S/2}$ are equilibria and $g\in H$,
from the elliptic regularity we have that $u_s, u_{s+S/2}\in H^4$.
On the other hand, we have
\begin{equation*}
A^2(u_s-u_{s+S/2}) = -A(f(u_s)-f(u_{s+S/2})).
\end{equation*}
Thus, recalling that $f\in C^3(\RR;\RR)$, we recover
\begin{equation}
\label{zelik5}
\|u_s-u_{s+S/2}\|_{H^5} \le C.
\end{equation}
On the other hand, by definition of $u_s$ and $u_{s+S/2}$, we have, for $\beta\in(0,1]$,
\begin{equation}
\|u(s+S/2)-u_s\|_{H^{1-\beta}} +\|u(s+S/2)-u_{s+S/2}\|_{H^{1-\beta}}
\le 2\sigma.
\end{equation}
This yields
\begin{equation}
\label{equi}
\|u_s-u_{s+S/2}\|_{H^{1-\beta}} \le 2\sigma,
\end{equation}
which, together with \eqref{zelik5} and interpolation, entails
\eqref{zelik4}.

We are now ready to construct the desired function
$\tilde u(t)$. Fix $S=2$ and introduce a cut-off function
$\theta\in C_0^\infty({\mathbb R})$ such that $\theta(t)\equiv0$ for
$t\le0$, $\theta(t)\equiv1$, for $t\ge1$, $0\le\theta(t)\le1$.
Then, for any $\sigma>0$ and any $N\in\mathbb{N}$, define a function $\tilde u(t)$ on the
interval $t\in[T-N,T-N+1]$, $T=T(\sigma)$, by the following formula
\begin{equation}
\label{zelik6}
\tilde u(t):=\theta(t-T+N)u_{T-N+1}+(1-\theta(t-T+N))u_{T-N}.
\end{equation}
This function is clearly smooth with respect to $t$ and
fulfills \eqref{Z4} and \eqref{Z5} (cf. \eqref{zelik4}).
Moreover, since both equilibria $u_{T-N+1}$ and
$u_{T-N}$ are close to $u(t)$ on the interval $t\in[T-N,T-N+1]$
(see \eqref{zelik3}), we also have, for all $t\le T$,
\begin{equation}
\label{zelik7}
\|u(t)-\tilde u(t)\|_{H^{1-\beta}}\le 2\sigma.
\end{equation}
Thus, on account of \eqref{zelik1} and \eqref{equi}, we infer that \eqref{Z1} holds.

To conclude the proof, it remains to show \eqref{Z2} and \eqref{Z3}.
By \eqref{Z5}, this is equivalent to check that the function
\begin{equation}
\label{zelik8}
\tilde \phi:=A^2\tilde u(t)+Af(\tilde u)-g
\end{equation}
is uniformly small. To this end, observe that
\begin{align}
\label{zelik9}
\tilde \phi(t)&=A^2u_{T-N}+\theta(t-T+N)A^2[u_{T-N+1}-u_{T-N}]\\
\nonumber
&+Af(u_{T-N}+\theta(t-T+N)[u_{T-N+1}-u_{T-N}])-g\\
\nonumber
&=\theta(t-T+N)A^2[u_{T-N+1}-u_{T-N}]\\
\nonumber
&+A[f(u_{T-N}+\theta(t-T+N)[u_{T-N+1}-u_{T-N}])-f(u_{T-N})]
\end{align}
and, using \eqref{zelik4} and the $C^3$-regularity of $f$,
it is not difficult to conclude that
\begin{equation}
\label{zelik10}
\| \tilde\phi(t)\|_H + \| \tilde\phi_t(t)\|_H \le C\sigma^\kappa.
\end{equation}
This inequality and \eqref{Z5} imply \eqref{Z2} and \eqref{Z3}.
\end{proof}

Now, for $L>0$ large enough (to be chosen below),
we look for a solution $v$ to the equation
\begin{equation}
\label{CH1-aux}
\epsi v_{tt} + v_{t} +A^2 v + A f(v) + LA^{-1} v = h(t),
\end{equation}
where
\begin{equation}
\label{h-def}
h(t)= g + LA^{-1} u(t).
\end{equation}
We recall that $U=(u,u_t)\in \mathcal{K}_\eb$ is given.

Observe that (cf. Corollary \ref{dissipation})
\begin{align*}
& \Vert h(T) \Vert_H^2 + \int_T^{T+1} \Vert A^{(1-\beta)/2}h_t(t) \Vert_H^2 dt \leq C_L (1+\Vert g\Vert_H^2),\\
& h_t \in C_b({\mathbb R};D(A^{(1-\beta)/2})),\quad \lim_{t\to -\infty}
\Vert h_t(t) \Vert_{D(A^{(1-\beta)/2})} =0,
\end{align*}
for any $\beta \in (0,1]$.

Let us prove the following (compare with \cite[Lemma 2.1]{Ze04})
\bete
\label{Zelik1}
Let $U=(u,u_t)\in \mathcal{K}_\eb$ be given. Then, for a sufficiently large $L>0$,
there exists a time $T=T(u,\epsi,L)<0$ such that equation
\eqref{CH1-aux}
possesses a unique strong backward solution $V \in L^\infty(-\infty,T;\calV_3)$ satisfying
\begin{equation}
\label{bound0}
\sqrt{\eps}\Vert v_{tt}(t)\Vert_H + \Vert v_t(t)\Vert_{H^2} + \Vert  v(t)\Vert_{H^4} \leq Q_L(\Vert g \Vert_H),
\qquad \forall\, t\le T,
\end{equation}
where the positive monotone function $Q_L(\cdot)$ is independent of $\epsi$.
Moreover, we have
\begin{equation}
\label{conv0}
\lim_{t\to -\infty} \Vert v_t(t)\Vert_{L^\infty(\Omega)}=0.
\end{equation}
\ente

\begin{proof}
We proceed as in \cite[proof of Lemma 2.1]{Ze04} with some modifications introduced by
Theorem \ref{Zelik}. We let $\sigma>0$ (possibly small enough) and look for a $v$ of the form
\begin{equation}
\label{vform}
v(t)=\tilde u(t) + w(t),
\end{equation}
where we recall that $\tilde u$ depends on $\sigma$.

Consequently, $w$ must solve the equation
\begin{equation}
\label{eqw}
\varepsilon w_{tt}+ w_t +A^2 w+A(f(\tilde u+w)-f(\tilde u)) + LA^{-1}w=\tilde h(t),
\end{equation}
where
\begin{equation}
\tilde h(t)= LA^{-1}(u - \tilde u)(t) - \phi(t).
\end{equation}
We will solve equation \eqref{eqw} by means of the inverse function Theorem
(cf., e.g., \cite[Thm.~2.1.2]{AP}).

Thus, we take $T<0$ small enough and consider the Banach space
$\Psi_b:=H^1_{loc}((-\infty,T];H)$  which is endowed with
the uniformly local norm
\begin{equation}\label{normG}
  \Vert G\Vert_{H^1_b((-\infty,T];H)}:=
   \sup_{t\in (-\infty,T-1]}\Vert G\Vert_{H^1((t,t+1);H)}.
\end{equation}
Then, it is easy to check that $\Psi_b$
is continuously embedded into $C_b((-\infty,T];H)$. We also define
\begin{equation}\label{smoothspace}
   \Phi_b := \big\{w\in C^2_b((-\infty,T];H) \cap C^1_b((-\infty,T];H^2) \cap C_b((-\infty,T];H^4):~
     \epsi w_{tt}+A^2w\in \Psi_b\big\}.
\end{equation}
Also the space $\Phi_b$ will be
endowed with the natural (graph) norm. Thus, we can define
the operator $\mathcal{T}:\Phi_b\to \Psi_b$
given by
\begin{equation*}
  \mathcal{T}_L: w \mapsto \epsi w_{tt}+ w_t +A^2 w+A(f(\tilde u+w)-f(\tilde u)) + LA^{-1}w.
\end{equation*}
Then, we observe that, recalling \eqref{Z1} and \eqref{Z3},
there exists $T=T(\sigma,u)\in {\mathbb R}$ such that
\begin{equation}
\label{convtrasl}
\Vert \tilde h\Vert_{H^1((t,t+1);H)} \leq C_L(\sigma + \sigma^\kappa), \qquad\forall\, t\leq T.
\end{equation}
Thus we work in a neighborhood of $0$ and we consider the linear operator $\mathcal{T}_L'(0)$
from $\Phi_b$ to $\Psi_b$ given by
\begin{equation*}
\mathcal{T}_L'(0)W =\epsi W_{tt}+ W_t +A^2 W+A(f'(\tilde u)W) + LA^{-1}W.
\end{equation*}
It then suffices to show that $\mathcal{T}_L'(0)$ is invertible. Thus we consider
the variation equation at $w=0$, namely,
\begin{equation}
\label{eqV}
\epsi W_{tt}+ W_t +A^2 W+A(f'(\tilde u)W) + LA^{-1}W=G(t).
\end{equation}
Here, we will assume that $G$ is a given function of $\Psi_b$.
We will now prove that \eqref{eqV} has a {\sl unique}\/ solution
$W\in \Phi_b$, that is, $\mathcal{T}_L'(0)W=G$,
for $T$ small enough and $L$ large enough.

Let us proceed formally by multiplying equation
\eqref{eqV} by $A^{-1}(W_t + \alpha W)$,
$\alpha>0$. We get
\begin{align}
\label{energy1}
& \frac{d}{dt} E_W
+ 2(1-\alpha\epsi)\Vert A^{-1/2}W_t\Vert_H^2 + 2\alpha\Vert A^{1/2}W\Vert_H^2
+ 2\alpha L\Vert A^{-1}W\Vert_H^2 + 2\alpha(f'(\tilde u)W,W)\\
\nonumber
&=  2(A^{-1/2}G,A^{-1/2}(W_t + \alpha W)) + (f^{\prime\prime}(\tilde u)\tilde u_t,W^2),
\end{align}
where
\begin{align}
E_W &= \eps\Vert A^{-1/2}W_t\Vert_H^2 + \Vert A^{1/2}W_H\Vert^2
+ L\Vert A^{-1}W\Vert_H^2 \\
\nonumber
&+ (f'(\tilde u)W,W) + 2\alpha\epsi(A^{-1/2}W,A^{-1/2}W_t)
+ \alpha\Vert A^{-1/2}W\Vert_H^2.
\end{align}

Observe that, recalling \eqref{f2} and using interpolation and Young's
inequality, we have
\begin{align}
\label{interp1}
&-(f'(\tilde u)W,W) \leq
\lambda\Vert W\Vert_H^2 \leq c_1\lambda\Vert W \Vert^{4/3}_{H^1}
\Vert W \Vert^{2/3}_{H^{-2}}\\
&\leq c_2\lambda\Vert A^{1/2} W \Vert_H^{4/3}\Vert A^{-1}W \Vert_H^{2/3}
\leq \frac{1}{2}\left(\Vert A^{1/2} W \Vert_H^{2}+
c_3\lambda^3\Vert A^{-1}W \Vert_H^{2}\right),
\nonumber
\end{align}
where $c_j$, $j=1,2,3$, are positive constants independent of $W$.

Then, recalling \eqref{f2}, we can choose $L>c_3\lambda^3$  so that
\begin{equation}
\label{interp2}
(f'(\tilde u)W,W) + \frac{1}{2}\left(\Vert A^{1/2}W\Vert_H^2 + L\Vert A^{-1}W\Vert_H^2\right)
\geq 0.
\end{equation}

Picking $\alpha$ small enough (but independent of $\epsi$ and $L$), we then deduce
\begin{equation}
C_1^{-1} \left(\epsi\Vert A^{-1/2}W_t\Vert_H^2 + \Vert A^{1/2}W\Vert_H^2\right)
\leq E_W \leq
C_1\left(\epsi\Vert A^{-1/2}W_t\Vert_H^2 + \Vert A^{1/2}W\Vert_H^2\right),
\end{equation}
for some $C_1>1$.

Consequently, from \eqref{energy1} we infer the inequality
\begin{equation}
\label{enineq1}
\frac{d}{dt} E_W + k E_W \leq C\Vert A^{-1/2}G\Vert_H^2
+ (f^{\prime\prime}(\tilde u)\tilde u_t,W^2) - \frac{\alpha}{2}\Vert A^{1/2}W\Vert_H^2,
\end{equation}
for some $k>0$.
On the other hand, for $\sigma$ (and $T$) small enough, we have, due to \eqref{Z5},
\begin{equation}
(f^{\prime\prime}(\tilde u)\tilde u_t,W^2) \leq \frac{\alpha}{2}\Vert A^{1/2}W\Vert_H^2.
\end{equation}
Thus we have
\begin{equation*}
\label{enineq1bis}
\frac{d}{dt} E_W + k E_W \leq C\Vert A^{-1/2}G\Vert_H^2.
\end{equation*}
Consider now a sequence $\{t_n\}\subset (-\infty,t)$ such that $t_n\to -\infty$.
Then Gronwall's inequality gives
\begin{align*}
\label{bound0}
\epsi\Vert A^{-1/2}W_t(t)\Vert_H^2 + \Vert A^{1/2}W(t)\Vert_H^2 &\leq
\left(\epsi\Vert A^{-1/2}W_t(t_n)\Vert_H^2 + \Vert A^{1/2}W(t_n)\Vert_H^2\right)e^{k(t_n-t)}\\
&+ C_L\int_{-\infty}^t e^{-k(t-s)}\Vert A^{-1/2}G(s)\Vert_H^2 ds.
\nonumber
\end{align*}
Since we are looking for solutions in $\Phi_b$ we can let $t_n$ go to $-\infty$ and recover
\begin{equation}
\label{bound1}
\epsi\Vert A^{-1/2}W_t(t)\Vert_H^2 + \Vert A^{1/2}W(t)\Vert_H^2 \leq
C_L\int_{-\infty}^t e^{-k(t-s)}\Vert A^{-1/2}G(s)\Vert_H^2 ds,
\end{equation}
for all $t\leq T$. Hence, in particular, $W$ is necessarily unique.

Let us now set
\begin{equation}
\widetilde W = W_t
\end{equation}
and observe that time differentiation of \eqref{eqV} gives
\begin{equation}
\label{eqVt}
\epsi\widetilde W_{tt}+ \widetilde W_t +A^2 \widetilde W+A(f'(\tilde u)\widetilde W)
+ LA^{-1}\widetilde W=G_t- A(f^{\prime\prime}(\tilde u)\tilde u_{t}W).
\end{equation}
Multiplying this equation by $A^{-1}(\widetilde W_t + \alpha \widetilde W)$ and arguing
as above, we find
\begin{align}
\label{bound2}
&\epsi\Vert A^{-1/2}\widetilde W_t(t)\Vert_H^2 + \Vert A^{1/2}\widetilde W(t)\Vert_H^2 \\
\nonumber
&\leq C_L\int_{-\infty}^t e^{-k(t-s)}
\left(\Vert A^{-1/2}G_t(s)\Vert_H^2
+ \Vert A^{1/2}(f^{\prime\prime}(\tilde u(s))\tilde u_{t}(s)W(s))\Vert_H^2\right) ds,
\end{align}
for all $t\leq T$. On account of \eqref{Z5} and \eqref{bound1}, observe now that
\begin{align}
\label{bound3}
 \Vert A^{1/2}(f^{\prime\prime}(\tilde u(s))\tilde u_{t}(s)W(s))\Vert_H
  & 
   \leq C \Vert A^{1/2}W(s)\Vert_H\\
\nonumber
  &\leq C_L
\sup_{t\in (-\infty,T-1]}\Vert G\Vert_{L^2((t,t+1);D(A^{-1/2}))}.
\end{align}
Therefore, \eqref{bound2} and \eqref{bound3} yield
\begin{equation}
\label{bound4}
\epsi\Vert A^{-1/2}W_{tt}(t)\Vert_H^2 + \Vert A^{1/2}W_t(t)\Vert_H^2 \leq
C_L \sup_{t\in (-\infty,T-1]}\Vert G\Vert^2_{H^1((t,t+1);D(A^{-1/2}))},
\end{equation}
and by comparison in \eqref{eqV}, we also deduce
\begin{equation}
\label{bound5}
\Vert A^{3/2}W(t)\Vert_H \leq
C_L \sup_{t\in (-\infty,T-1]}\Vert G\Vert_{H^1((t,t+1);D(A^{-1/2}))}.
\end{equation}

We can now multiply equation \eqref{eqVt} by $\widetilde W_t + \alpha \widetilde W$
and use the identity
\begin{align*}
 (A(f^{\prime}(\tilde u) \widetilde W), \widetilde W_t)&=\int_\Omega (\nabla(f^{\prime}(\tilde u) \widetilde W),
      \nabla \widetilde W_t)\\
 &= \frac{d}{dt}\int_\Omega f^{\prime}(\tilde u)  \frac{\vert \nabla \widetilde W \vert^2}{2}
   - \int_\Omega f^{\prime\prime}(\tilde u)\tilde u_t \frac{\vert \nabla \widetilde W \vert^2}{2}\\
 &-\int_\Omega \nabla\cdot[f^{\prime\prime}(\tilde u)\widetilde W\nabla\tilde u]\widetilde W_t.
\end{align*}
This gives the identity
\begin{align}
\label{energy2}
& \frac{d}{dt} \left[\eps\Vert \widetilde W_t\Vert_H^2 + \Vert A\widetilde W\Vert_H^2
+ L\Vert A^{-1/2}\widetilde W\Vert_H^2
+ 2\alpha\epsi(\widetilde W,\widetilde W_t)
+ \alpha\Vert \widetilde W\Vert_H^2 + \int_\Omega f^{\prime}(\tilde u) \vert \nabla \widetilde W \vert^2\right]\\
\nonumber
& + 2(1-\alpha\epsi)\Vert \widetilde W_t\Vert_H^2 + 2\alpha\Vert A \widetilde W\Vert_H^2
+ 2\alpha L\Vert A^{-1/2}\widetilde W\Vert_H^2
+ 2\alpha(A^{1/2}(f'(\tilde u)\widetilde W),A^{1/2}\widetilde W)\\
\nonumber
&=  2(G_t,\widetilde W_t + \alpha \widetilde W)
+ \int_\Omega f^{\prime\prime}(\tilde u)\tilde u_t  \vert \nabla \widetilde W \vert^2
+ 2\int_\Omega \nabla\cdot[f^{\prime\prime}(\tilde u)\widetilde W\nabla\tilde u]\widetilde W_t
- 2(A(f^{\prime\prime}(\tilde u)\tilde u_t W), \widetilde W_t + \alpha \widetilde W).
\end{align}
Observe that, by \eqref{f2} and \eqref{Z4}-\eqref{Z5},
\begin{align}
\label{ineq1bis}
&\int_\Omega f^{\prime}(\tilde u) \vert \nabla \widetilde W \vert^2
\leq C\Vert A^{1/2}\widetilde  W\Vert_H^2,\\
\label{ineq1}
&(A^{1/2}(f^{\prime}(\tilde u)\widetilde W), A^{1/2}\widetilde W)
\leq C\Vert A^{1/2}\widetilde  W\Vert_H^2,\\
\label{ineq1quater}
&\int_\Omega f^{\prime\prime}(\tilde u)\tilde u_t  \vert \nabla \widetilde W \vert^2
\leq C\Vert A^{1/2}\widetilde  W\Vert_H^2,\\
\label{ineq1ter}
&\int_\Omega \nabla\cdot[f^{\prime\prime}(\tilde u)\widetilde W\nabla\tilde u]\widetilde W_t
\leq C\Vert A^{1/2} \widetilde W\Vert_H \Vert \widetilde W_t \Vert_H,\\
\label{ineq2}
&\Vert A(f^{\prime\prime}(\tilde u)\tilde u_t W)\Vert_H \leq C\Vert AW \Vert_H.
\end{align}
Next, due to
\begin{equation*}
\Vert A^{1/2}\widetilde W\Vert^2 \leq C \Vert \widetilde W\Vert^{4/3}_{H^2}
\Vert \widetilde W \Vert^{2/3}_{H^{-1}}
\leq C \Vert A \widetilde W \Vert_H^{4/3}\Vert A^{-1/2}\widetilde W \Vert_H^{2/3},
\end{equation*}
also in this case we can choose $L$ large enough so that
\begin{equation*}
\int_\Omega f^{\prime}(\tilde u) \vert \nabla \widetilde W \vert^2
+ \frac{1}{2}\left(\Vert A \widetilde W\Vert_H^2 + L\Vert A^{-1/2}\widetilde W\Vert_H^2\right)
\geq 0.
\end{equation*}
We now indicate by $\widetilde E_{\widetilde W}$ the expression within square brackets
in identity \eqref{energy2}.
Then we can choose $\alpha$ small enough (but independent of $\epsi$ and $L$)
in such a way that
\begin{equation}\label{newlabel}
C_2^{-1} \left(\epsi\Vert \widetilde W_t\Vert_H^2 + \Vert A\widetilde W\Vert_H^2\right)
\leq \widetilde E_{\widetilde W} \leq
C_2\left(\epsi\Vert \widetilde W_t\Vert_H^2 + \Vert A\widetilde W\Vert_H^2\right),
\end{equation}
for some $C_2>1$. Then, from  \eqref{energy2} and \eqref{ineq1}-\eqref{ineq2} we deduce the inequality
\begin{equation}
\label{enineq2}
\frac{d}{dt} \widetilde E_{\widetilde W} + \alpha \widetilde E_{\widetilde W}
\leq C(\Vert G_t\Vert_H^2 + \Vert AW \Vert_H^2 + \Vert A^{1/2}W_t\Vert_H).
\end{equation}
Therefore, owing to \eqref{newlabel} and recalling \eqref{bound4},
a further application of Gronwall's inequality to \eqref{enineq2} yields
\begin{equation}
\label{bound6}
\epsi\Vert \widetilde W_t(t)\Vert_H^2 + \Vert A \widetilde W(t)\Vert_H^2
\leq C_L\int_{-\infty}^t e^{-\alpha(t-s)}(\Vert G_t(s)\Vert_H^2 + \Vert AW(s) \Vert_H^2) ds,
\end{equation}
for all $t\leq T$. This yields, on account of \eqref{bound4} and \eqref{bound5},
\begin{equation}
\label{bound7}
\epsi\Vert W_{tt}(t)\Vert_H^2 + \Vert A W_t(t)\Vert_H^2
\leq C_L \Vert G\Vert^2_{H^1_b((-\infty,T);H)},
\end{equation}
and suitable comparison arguments in \eqref{eqV} entail that
\begin{equation}\label{bound8}
  \Vert W(t)\Vert_{H^4} \leq C_L \Vert G\Vert_{H^1_b((-\infty,T);H)}
\end{equation}
and, correspondingly,
\begin{equation}
\label{bound8bis}
  \Vert W_{tt}+A^2 W(t)\Vert_{H^1_b((-\infty,T);H)}
   \leq C_L \Vert G\Vert_{H^1_b((-\infty,T);H)}.
\end{equation}
A priori estimates \eqref{bound7} and \eqref{bound8}-\eqref{bound8bis}
combined with standard arguments allow
us to conclude that equation \eqref{eqV} has indeed a (unique) solution in $\Phi_b$.

Thus, owing to the inverse function theorem, we conclude that \eqref{eqw}
has a unique solution $w$ for $T$ small enough such that
\begin{equation}
\label{bound9}
\Vert  w_t(t)\Vert_{H^2} + \Vert  w(t)\Vert_{H^4}
\leq C_L (\Vert g \Vert_H) .
\end{equation}
where $C_L$ is independent of $\epsi$ (cf. \eqref{bound7}).

Moreover, thanks to \eqref{convtrasl}, if we fix a
neighborhood ${\mathcal W}$ of $0$  in the Banach space \eqref{smoothspace},
then we can choose $T$ small enough such that $(w(t),w_t(t))\in {\mathcal W}$
for all $t\leq T$. Thus, recalling \eqref{vform} and Theorem \ref{Zelik}, we have, in particular,
\begin{equation}
\label{conv1}
\lim_{t \to -\infty}\Vert v_t(t)\Vert_{H^2} =0
\end{equation}
and can conclude  that we have found a backward strong solution
to \eqref{CH1-aux} which satisfies \eqref{bound0} and, on account of the arbitrariness of $\sigma$ (cf. also \eqref{Z5}),
\eqref{conv0}.
\end{proof}

We are now in a position to prove the main result of this section, namely,
\bete
\label{backsmooth}
Let $U=(u,u_t)\in {\mathcal K}_\eb$. Then there exists $T=T_u$ such that
\begin{equation}
\label{smooth}
U\in C_b((-\infty,T],\calV_3),
\end{equation}
and, for all $t\leq T_u$,
\begin{equation}
\label{backest}
\Vert u(t)\Vert_{H^4}^2+\|u_t(t)\|_{H^2}^2+\eb\|u_{tt}(t)\|^2_{H} \leq Q(\Vert g\Vert_H).
\end{equation}
\ente

\begin{proof}
The proof follows closely \cite{Ze04} with suitable adaptations.
Our goal is to prove that
\begin{equation}
\label{basicid}
U(t)\equiv V(t),
\end{equation}
for all $t\leq T$ and some $T<0$, where $V$ is given by Theorem~\ref{Zelik1}. According to Theorem \ref{Zelik0}, we consider
a sequence $\{U^{n_k}(t)\}$ for $t\geq t_k$ satisfying \eqref{CH1-nk}-\eqref{iniz-nk}
and such that \eqref{bdd-conv-nk} holds. Also, we define
\begin{equation}
v^{n_k}(t) = P_{n_k}v(t), \qquad \forall\,t\leq T,
\end{equation}
where $T$ is given by Theorem~\ref{Zelik1}, and we observe that, due to Theorem \ref{Zelik1}, we have
\begin{equation}
\label{convPr}
\lim_{k\to\infty}\Vert V^{n_k} - V\Vert_{C_b((-\infty,T];\calV_2)} = 0, \quad
\lim_{k\to\infty}\Vert v^{n_k} - v\Vert_{C_b((-\infty,T]: C^1(\bar\Omega))} = 0.
\end{equation}
In addition, we have (cf. \eqref{bound0})
\begin{equation}
\label{convPrt}
\lim_{k\to\infty}\Vert v^{n_k}_t - v_t\Vert_{C_b((-\infty,T]\times\bar\Omega)} = 0.
\end{equation}
Let us now set
\begin{equation}
\label{Zdef}
Z(t)= u(t) -v(t), \quad Z^{n_k}(t) = u^{n_k}(t) -v^{n_k}(t),
\end{equation}
and observe that
\begin{align}
\label{eqZ}
&\disp \epsi Z^{n_k}_{tt} + Z^{n_k}_{t} +A^2 Z^{n_k}
+ A P_{n_k}(f(v^{n_k} + Z^{n_k}) - f(v^{n_k}))
+ LA^{-1}Z^{n_k}  = h^{n_k} , \\
\label{eqZin}
&\disp Z^{n_k}|_{t=t_k}=u^k_{0} - P_{n_k}v(t_k), \quad
Z^{n_k}_t|_{t=t_k}=u^k_{1}-P_{n_k}v_t(t_k),
\end{align}
where
\begin{equation}
h^{n_k}(t) := A P_{n_k}(f(v(t)) - f(v^{n_k}(t))).
\end{equation}
Observe also that, on account of \eqref{convPr}, there holds
\begin{equation}
\label{convh}
\lim_{k\to\infty}\Vert A^{-1/2}h^{n_k}\Vert_{C_b((-\infty,T]\times\bar\Omega)} =0.
\end{equation}
Moreover, thanks to Theorem~\ref{Zelik0} and \eqref{bound0}, we also have
\begin{equation}
\label{bddini}
\Vert (Z^{n_k}(t_k),Z^{n_k}_t(t_k)) \Vert_{\calV_0} \leq C, \qquad\forall\,k\in\mathbb{N}.
\end{equation}
If we multiply equation \eqref{eqZ} by $A^{-1}(Z^{n_k}_t + \alpha Z^{n_k})$, then we obtain
\begin{equation}
\label{energy3}
\frac{d}{dt} E_{Z^{n_k}} + \alpha E_{Z^{n_k}} = H_{n_k}, \qquad\textrm{ in } (-\infty,T],
\end{equation}
where
\begin{align}
\label{energyk}
E_{Z^{n_k}} &=  \eps\Vert A^{-1/2}Z^{n_k}_t\Vert_H^2 + \Vert A^{1/2} Z^{n_k}\Vert_H^2
+ L\Vert A^{-1} Z^{n_k}\Vert_H^2 + 2\alpha\epsi( Z^{n_k}, Z^{n_k}_t) \\
\nonumber
&+ \alpha\Vert A^{-1/2}Z^{n_k}\Vert_H^2
+ 2\left(F(v^{n_k} + Z^{n_k}) - F(v^{n_k}) - f(v^{n_k})Z^{n_k},1\right)
\end{align}
and
\begin{align}
\label{sourcek}
H_{n_k} &:=  -(2-3\eps)\Vert A^{-1/2} Z^{n_k}_t\Vert_H^2 -\alpha \Vert A^{1/2} Z^{n_k}\Vert_H^2
- \alpha L\Vert A^{-1} Z^{n_k}\Vert_H^2\\
\nonumber
&+2\alpha\left(F(v^{n_k} + Z^{n_k}) - F(v^{n_k}) - f(v^{n_k})Z^{n_k}
- (f(v^{n_k}+ Z^{n_k}) - f(v^{n_k}))Z^{n_k},1\right)\\
\nonumber
&+ 2\alpha^2\epsi( Z^{n_k}, Z^{n_k}_t) + \alpha^2\Vert Z^{n_k}\Vert_H^2
+ 2(A^{-1/2}h^{n_k},A^{-1/2}(Z^{n_k}_t + \alpha Z^{n_k}))\\
\nonumber
&+ 2\left(f(v^{n_k} + Z^{n_k}) - f(v^{n_k}) - f'(v^{n_k})Z^{n_k},v^{n_k}_t\right).
\end{align}
We can now take advantage of \cite[(2.52)-(2.53)]{Ze04} to estimate
the above nonlinear terms. Then, also recalling the choice of $L$ and
\eqref{interp2}, we can find some positive constants $\alpha_1$,
$C_1$ and $C_2$ (all independent of $v^{n_k}$, $Z^{n_k}$, $k$,
$L$ and $\eps$) such that
\begin{align}
\label{sourcekest}
H_{n_k} &\leq -\frac{\alpha_1}{2}
\left(\Vert A^{1/2} Z^{n_k}\Vert_H^2 + \Vert A^{-1/2} Z^{n_k}_t\Vert_H^2\right)
- \frac{\alpha_1 L}{2} \Vert A^{-1} Z^{n_k}\Vert_H^2 -2\alpha_1(\vert Z^{n_k}\vert^{p+4},1)\\
\nonumber
&+ C_1\Vert A^{-1/2} h^{n_k}\Vert_H^2 + C_2\Vert v^{n_k}_t \Vert_{L^\infty(\Omega)}
\left(\vert Z^{n_k}\vert^2(1 + \vert v^{n_k}\vert^{p+1} + \vert Z^{n_k}\vert^{p+1}),1\right).
\end{align}
Using \eqref{conv0} and \eqref{convPr}-\eqref{convPrt}, we can also find $T'\leq T$
such that
\begin{equation}
H_{n_k}(t) \leq  C_1\Vert A^{-1/2} h^{n_k}(t)\Vert_H^2, \qquad \forall\,t\leq T'.
\end{equation}
Then, applying Gronwall's lemma to \eqref{energy3} and using \cite[(2.51)]{Ze04}
and \eqref{interp1}, we get
\begin{align}
&\eps\Vert A^{-1/2}Z^{n_k}_t(t)\Vert_H^2 + \Vert A^{1/2} Z^{n_k}(t)\Vert_H^2\\
\nonumber
&\leq C_3\left(1+ \Vert (Z^{n_k}(t_k),Z^{n_k}_t(t_k)) \Vert^2_{\calV_0}\right)e^{-\alpha(t-t_k)}
+2C_1 \int_{t_k}^t e^{-\alpha(t-s)}\Vert A^{-1/2} h^{n_k}(s)\Vert_H^2 ds,
\end{align}
for all $t\leq T'$, where $C_3$ is also independent of $k$. 
Then, on account of \eqref{convh} and \eqref{bddini}, we
let $k$ go to $\infty$ and we recover (cf. \eqref{Zdef})
\begin{equation}
\eps\Vert A^{-1/2}Z_t(t)\Vert_H^2 + \Vert A^{1/2} Z(t)\Vert_H^2 \leq 0, \qquad\forall\,t\leq T'.
\end{equation}
Thus \eqref{basicid} is proven. Estimate \eqref{backest} follows from \eqref{bound0} and
\eqref{basicid}. The proof is complete.
\end{proof}

We conclude this section by proving that the solution $U=(u,u_t)\in {\mathcal K}$ is unique
until it is regular (cf. \cite[Thm.~2.2]{Ze04}).

\bete
\label{unismooth}
Let $U=(u,u_t)\in {\mathcal K}_\eb$ satisfy \eqref{backest} for every $t\leq T_u$. Consider another
complete energy solution $\tilde U=(\tilde u,\tilde u_t)\in {\mathcal K}_\eb$ such that, for some $\tilde T<T_u$ and
for all $t\leq \tilde T <T_u$,
\begin{equation}
\label{unism1}
U(t)\equiv \tilde U(t).
\end{equation}
Then we necessarily have
\begin{equation}
\label{unism2}
U(t)\equiv \tilde U(t),
\end{equation}
for all $t\leq T_u$.
\ente

\begin{proof}
Arguing as in the  proof of Theorem \ref{backsmooth}, we consider a sequence of Galerkin solutions
$\{\tilde U^{n_k}(t)\}$ which approximate $\tilde U$ for $t\geq t_k$. Namely,
we assume \eqref{CH1-nk}, \eqref{iniz-nk} and \eqref{bdd-conv-nk}
to be satisfied. Then, we set
\begin{equation*}
U^{n_k} = P_{n_k}U, \quad z=u-\tilde u,
\quad z^{n_k}=u^{n_k}-\tilde u^{n_k},
\end{equation*}
where $u^{n_k}= P_{n_k} u$.
Thus, similarly to \eqref{eqZ}-\eqref{eqZin}, we have
\begin{align}
\label{eqZ2}
&\disp \epsi z^{n_k}_{tt} + z^{n_k}_{t} +A^2 z^{n_k}
+ A P_{n_k}(f(\tilde u^{n_k} + z^{n_k}) - f(\tilde u^{n_k}))
+ LA^{-1}z^{n_k}  = h^{n_k}(t), \\
\label{eqZ2in}
&\disp z^{n_k}|_{t=t_k}=u^k_{0} - P_{n_k}\tilde u(t_k), \quad
z^{n_k}_t|_{t=t_k}=u^k_{1}-P_{n_k}\tilde u_t(t_k),
\end{align}
where
\begin{equation}\label{hhhh}
h^{n_k}(t) := A P_{n_k}(f(\tilde u(t)) - f(\tilde u^{n_k}(t))) + LA^{-1}z^{n_k}.
\end{equation}
for some positive $L$ which will be chosen below.

Observe that (cf. \eqref{convh}),
\begin{equation}
\label{convh1}
\lim_{k\to\infty}\Vert A^{-1/2}(h^{n_k}-LA^{-1}z^{n_k})\Vert_{C_b((-\infty,T_u];H)} =0,
\quad \Vert A^{-1/2}h^{n_k}\Vert_{C_b((-\infty,T_u];H)}\leq C_4,
\end{equation}
for some $C_4>0$ independent of $k$.

If we multiply equation \eqref{eqZ2} by $A^{-1}(z^{n_k}_t + \alpha z^{n_k})$, then we obtain
(cf. \eqref{energy3})
\begin{equation}
\label{energy3bis}
\frac{d}{dt} E_{z^{n_k}}(t) + \alpha E_{z^{n_k}}(t) = H_{n_k}(t), \qquad \forall\,t\leq T_u.
\end{equation}
Here $E_{z^{n_k}}$ and $H_{n_k}$ are defined as in \eqref{energyk}
and \eqref{sourcek}, respectively. Thus, recalling \eqref{sourcekest},
we can find some positive constants
$\alpha_2$, $C_5$ and $C_6$ (all independent of $\tilde u^{n_k}$, $z^{n_k}$,
$k$, $L$ and $\eps$) such that
\begin{align}
\label{sourcekest2}
H_{n_k} &\leq -\frac{\alpha_2}{2}
\left(\Vert A^{1/2} z^{n_k}\Vert_H^2 + \Vert A^{-1/2} z_t^{n_k}\Vert_H^2\right)
- \frac{\alpha_2 L}{2} \Vert A^{-1} z^{n_k}\Vert_H^2 -2\alpha_2(\vert z^{n_k}\vert^{p+4},1)\\
\nonumber
&+ C_5\Vert A^{-1/2} h^{n_k}\Vert_H^2 + C_6\Vert  u^{n_k}_t \Vert_{L^\infty(\Omega)}
\left(\vert z^{n_k}\vert^2(1 + \vert u^{n_k}\vert^{p+1} + \vert z^{n_k}\vert^{p+1}),1\right),
\end{align}
in $(-\infty,T_u]$. Using now \eqref{backest}, we get
\begin{align}
\label{sourcekest3}
H_{n_k} &\leq -\frac{\alpha_2}{2}
\left(\Vert A^{1/2} z^{n_k}\Vert_H^2 + \Vert A^{-1/2} z_t^{n_k}\Vert_H^2\right)
- \frac{\alpha_2 L}{2} \Vert A^{-1} z^{n_k}\Vert_H^2 -2\alpha_2(\vert z^{n_k}\vert^{p+4},1)\\
\nonumber
&+ C_5\Vert A^{-1/2} h^{n_k}\Vert_H^2 + C_7
\left(\vert z^{n_k}\vert^2(1 + \vert z^{n_k}\vert^{p+1}),1\right),\\
\nonumber
&\leq -\frac{\alpha_2}{2}
\left(\Vert A^{1/2} z^{n_k}\Vert_H^2 + \Vert A^{-1/2} z_t^{n_k}\Vert_H^2\right)
- \frac{\alpha_2 L}{2} \Vert A^{-1} z^{n_k}\Vert_H^2 -\alpha_2(\vert z^{n_k}\vert^{p+4},1)\\
\nonumber
&+ C_5\Vert A^{-1/2} h^{n_k}\Vert_H^2 + C_8\Vert z^{n_k}\Vert_H^2.
\end{align}
Using \eqref{interp1} and Young's inequality, we can then choose $L$ large enough so that
\begin{equation}
\label{sourcekest4}
H_{n_k}(t) \leq C_5\Vert A^{-1/2} h^{n_k}(t)\Vert_H^2,\qquad \forall\,t\leq T_u.
\end{equation}
Applying Gronwall's inequality to \eqref{energy3bis} we obtain,
thanks to \eqref{sourcekest4},
\begin{equation}
\label{energy3tris}
E_{z^{n_k}}(t) \leq E_{z^{n_k}}(t_k)e^{-\alpha(t-t_k)}
+ C_5\int_{t_k}^t e^{-\alpha(t-s)}\Vert A^{-1/2} h^{n_k}(s)\Vert_H^2  ds,
\qquad \forall\,t\leq T_u.
\end{equation}
We now let $k$ go $+\infty$ in \eqref{energy3tris} and recalling \eqref{unism1}
and \eqref{convh1}, we obtain
\begin{equation}
\label{energy4}
\Vert Z(t)\Vert_0 \leq
C_5L^2\int_{\tilde T}^t e^{-\alpha(t-s)}\Vert A^{-3/2} z(s)\Vert_H^2  ds,
\qquad \forall\,t\in [\tilde T, T_u].
\end{equation}
Finally, a further application of Gronwall's lemma to \eqref{energy4} entails
\eqref{unism2}.
\end{proof}


\section{Smoothness of global attractors}
\label{smoothattr}
We can now show that, under suitable conditions, the trajectory attractor
${\mathcal A}_\epsi^{tr}$ (see \eqref{trajattr}) consists of strong solutions.
Let us first analyze the two-dimensional case in the case of cubic controlled growth.
Namely, we take $f\in C^3(\RR;\RR)$ with $f(0)=0$ satisfying \eqref{f1}, \eqref{f2} and
\begin{equation}
\label{f3bis}
f'''\in L^\infty(\RR;\RR).
\end{equation}
We know that, for any $\epsi>0$, problem $P_\epsi$ generates a
(strongly continuous) semigroup $S_\epsi(t)$ on $\calV_2$ (see
\cite[Thm.~2.2]{GSZ}, cf. also \cite[Rem.~2.1]{GSZ}), i.e., the trajectories are quasi-strong solutions.
Moreover, the dynamical system
$(\calV_2,S_\epsi(t))$ possesses a global attractor
$\mathbb{A}_\epsi$ which is bounded in $\calV_3$ (cf. \cite[Thm.~4.1]{GSZ}). On the other hand,
in the case of energy solutions, we also know that $P_\epsi$
generates a (strongly continuous) semigroup $\tilde S_\epsi(t)$ on
$\calV_0$ which has the global attractor $\mathcal{A}_\eps$ (see
\cite[Sec.~6]{GSZ}). Thus, on account of Theorem \ref{Zelik0}, we
have that ${\mathcal K}$ coincides with the set of all the complete
energy bounded solutions and
\begin{equation}
\label{equalattr}
\mathcal{A}_\epsi\equiv \Pi_0{\mathcal A}^{tr}_\eb,
\end{equation}
where $\Pi_0U(t)=U(0)$.

It is clear that $\mathbb{A}_\epsi\subset {\mathcal A}_\epsi$, while the
validity of the opposite inclusion
was an open question so far. Indeed we are now ready to prove the following

\bete
\label{smooth2D}
Let \eqref{f1}, \eqref{f2}, \eqref{regog}  and \eqref{f3bis} hold.
Then
\begin{equation}
\label{smooth2Deq}
\mathbb{A}_\epsi\equiv {\mathcal A}_\epsi.
\end{equation}
\ente

\begin{proof}
Let us consider first the set ${\mathbb K}_\eps$ of all the complete
bounded strong solutions so that $\mathbb{A}_\eps\equiv
\Pi_0{\mathbb K}_\eps.$ It is clear that ${\mathbb K}_\eps
\subset {\mathcal K}_\eps $. Let us consider now a bounded
complete energy solution $U$ to \eqref{CH1-regular}. Then, by
Theorem \ref{backsmooth}, there exists a time $T=T_u$ such that
$U(t) \in \calV_3$ for all $t\leq T$ and satisfies a bound like
\eqref{backest}. On the other hand, we know that there exists a
unique strong solution $\tilde U$ to \eqref{CH1-regular} bounded on
$[T,+\infty)$ and such that $\tilde U(T) = U(T)$. Thus we can
construct a complete strong solution
\begin{equation}
\label{complsol}
U^*(t) = \begin{cases}
&\tilde U(t), \qquad t>T,\\
&U(t), \qquad t\leq T,
\end{cases}
\end{equation}
which is bounded on $\RR$ in the $\calV_3$-norm. Thus, Theorem
\ref{unismooth} implies that $U\equiv U^*$. Therefore, we have that
${\mathcal K}_\eps \subset {\mathbb K}_\eps $ and the proof is
complete.
\end{proof}

In the three-dimensional case or in two dimensions with supercubical growth,
we cannot say more than Theorem
\ref{Zelik0} about the existence of global solutions without making a
restriction on $\eps$. More precisely, we have the following result
(see \cite[Thm.~2.4]{GSSZ})

\bete\label{globexuni} Let \eqref{f1}-\eqref{f3} and
\eqref{regog} hold.
 Then, there exist  $\epsi_0>0$ and a nonincreasing
 positive function $R:(0,\epsi_0)\to (0,+\infty)$
 with the property
 \begin{equation}\label{propR}
   \lim_{\epsi\searrow 0}R(\epsi)= +\infty,
 \end{equation}
 such that, for every $\epsi\in(0,\epsi_0)$
 and every initial condition
 $U_0=(u_0,u_1)\in \mathcal{V}_1$
 satisfying
 \begin{equation}\label{cond-dati-piccoli}
   \|U_0\|_{\calV_1} \le R(\epsi),
 \end{equation}
 there exists a (unique) global weak solution $U=(u,u_t)$
 to problem $P_\epsi$ such that
 \begin{equation}\label{stimadiss}
   \|U(t)\|_{\calV_1}\le Q(\|U_0\|_{\calV_1}) e^{-\kappa t} +Q(\|g\|_H), \quad\perogni t\geq 0,
 \end{equation}
 for some $\kappa>0$ and some positive increasing
 monotone function $Q$ both independent of $\epsi$.
\ente

Following \cite{Ze04}, we consider the $\calV_1$-ball
$B(R(\eps))$ of radius $R(\eps)$ centered at $0$ and define
the solving operator
\begin{equation*}
S_\eb(t): B(R(\eps)) \to \calV_1, \qquad U(t)=S_\eb(t)U(0),
\quad \forall\,t\geq 0.
\end{equation*}
Thanks to estimate \eqref{stimadiss}, we have
\begin{equation*}
\Vert S_\eb(t)(B(R(\eps))) \Vert_{\calV_1} \leq  \Lambda,\quad\forall\,t\geq 0,
\end{equation*}
for an appropriate positive quantity $\Lambda$
(clearly also depending on $\|g\|_H$). Then, we set
\begin{equation}\label{inclB}
  \mathbb{B}_{\eps}:=
   \left[\bigcup_{t\geq 0} S_\eb(t)(B(R(\eps)))\right]_{\calV_1},
\end{equation}
where $[\cdot]_{\calV_1}$ denotes the closure in $\calV_1$, and we
note that $\mathbb{B}_{\eps}$ is bounded in $\calV_1$ by $\Lambda$.
Moreover, $\mathbb{B}_{\eps}$ is closed and positively
invariant by construction.
Thus, we can say that $(\mathbb{B}_{\eps},S_\eps(t))$
is a dissipative dynamical system.

On the other hand, we know that (see \cite[Thm.~2.7]{GSSZ}),
up to possibly taking a smaller $\epsi_0$ and correspondingly
modifying the function $R$, for all $\eps\in (0,\eps_0)$
$(\mathbb{B}_{\eps},S_\eps(t))$ possesses a global
attractor $\mathbb{A}_\eps$ bounded in $\calV_3$. Consider now
the set ${\mathbb K}_\eps$ of all complete and bounded weak
solutions taking values in the phase space $\mathbb{B}_{\eps}$.
Thanks to the above considerations, we have that
$\mathbb{A}_\eps=\Pi_t{\mathbb K}_\eps$, for any
$t\in\RR$. Hence, each $U\in {\mathbb K}_\eps$ is indeed a global
strong solution. We can now prove

\bete \label{smooth3D}
 Let the assumptions of
 Theorem~\ref{globexuni} hold. Then there exists $\eps_1\in
 (0,\eps_0]$ such that, if $\eps\in (0,\eps_1)$, then
 \begin{equation}\label{smooth3Deq}
  {\mathbb K}_\eps \equiv {\mathcal K}_\eb,
 \end{equation}
 where ${\mathcal K}_\eb$ is defined in Theorem \ref{Zelik0}.
 In particular, we have
 \begin{equation}\label{smooth3Deq2}
  \mathbb{A}_\epsi\equiv \Pi_0{\mathcal A}^{tr}_\eb.
 \end{equation}
\ente

\begin{proof} It is clear that
${\mathbb K}_\eps \subset {\mathcal K}_\eb$. On the other hand, if we
consider a bounded complete energy solution $U$ given by
Theorem~\ref{Zelik0} (i.e., an element of the trajectory attractor),
then Theorem \ref{backsmooth} implies the existence
of a time $T=T_u$ such that $U(t) \in \calV_3$ for all $t\leq T$ and
satisfies a bound like \eqref{backest}. In particular, this bound
entails that, for all $t\leq T_u$,
\begin{equation}
\label{backest2}
\Vert U(t)\Vert_{\calV_1} = \left(\|u(t)\|_{H^2}^2
   +\epsi\|u_t(t)\|_H^2\right)^{1/2} \leq Q_1(\Vert g\Vert_H),
\end{equation}
where $Q_1(\cdot)$ is a computable function
whose expression is independent of $\eps$. At this point, we
restrict ourselves to those $\eps\in (0,\eps_0)$ such that
\begin{equation}\label{radiusest}
  Q_1(\Vert g\Vert_H) \leq R(\eps).
\end{equation}
By \eqref{propR} and the (decreasing) monotonicity of $R(\cdot)$,
this will hold for all $\epsi$ in some
interval $(0,\epsi_1)$, where $\epsi_1\le\epsi_0$.
By Theorem~\ref{globexuni}, we know that there
exists a unique solution $\tilde U\in
C_b([T,+\infty);\calV_1)$ to \eqref{CH1-regular} such that $\tilde
U(T) = U(T)$. Thus, recalling \eqref{complsol},
we can construct once more a
complete weak solution $U^*$. Moreover,
being $B(R(\eps))\subset \mathbb{B}_{\eps}$
by \eqref{inclB}, this solution is such that
\begin{equation*}
  U^*(t) \in \mathbb{B}_{\eps} \quad\perogni
   t\in\RR.
\end{equation*}
Thus we have that necessarily $U^*\in {\mathbb K}_\eps $.
On the other hand, Theorem \ref{unismooth}
implies that $U\equiv U^*$. Therefore, we have once again
${\mathcal K}_\eps \subset {\mathbb K}_\eps $.
\end{proof}

\section{Exponential regularization of energy solutions}\label{s5.exp}

In this  section, we use the proved regularity of the trajectory
attractor $\mathcal A_\eb^{tr}$ in order to verify that every energy
solution is a sum of an exponentially decaying and a smooth function.
This fact, together with the transitivity of exponential attraction,
will allow us to establish that the exponential attractor $\mathcal
M_\eb$, constructed in  \cite{GSSZ} for the case of {\it weak}
solutions (and small $\eb>0$) is automatically the exponential
attractor in the class of {\it energy} solutions as well. Moreover, that
exponential regularization property will help us to verify the
well-posedness and dissipativity of the 2D problem \eqref{CHin} in
the intermediate (between $\mathcal V_0$ and $\mathcal V_2$)
phase space $\mathcal V_1$. To this end, we first formulate  the
regularity estimate for the solutions on the  attractor $\mathcal
A_\eb^{tr}$ obtained before.

\begin{coro}\label{Cor4.attr} Let the assumptions of
Theorem \ref{smooth2D} or Theorem \ref{smooth3D} be satisfied.
Then any complete trajectory $U\in\mathcal K_\eb$ satisfies
\begin{equation}\label{4.smattr}
\|u(t)\|_{H^4}^2+\|u_{t}(t)\|_{H^2}^2+\eb\|u_{tt}(t)\|_{H}^2\le C,
\end{equation}
where the constant $C$ is independent of $t$, $\eb$ and of the
concrete choice of the trajectory $u$.
\end{coro}
The next result, which gives the analogue of Theorem \ref{Zelik} for
the forward in time energy solutions, is the main technical tool of this
section.

\begin{lemm}\label{Lem4.appr} Let the assumptions of Corollary \ref{Cor4.attr} hold,
so that the trajectory attractor $\mathcal A_\eb^{tr}$ of equation
\eqref{CHin} possesses the regularity
 \eqref{4.smattr}. In addition, suppose that $f\in C^4(\RR;\RR)$.
 Then, for every energy solution $U\in K^+_\eb$ and
 every $\sigma>0$, there exist $T=T(\eb,\kappa, M^\eb_u(0))$
 and a function $\tilde u(t)$ such that
 \begin{equation}\label{4.sm}
 \|\tilde u(t)\|_{H^4}^2+\|\tilde u_{t}(t)\|_{H^2}^2
 +\eb\|\tilde u_{tt}(t)\|_{H}^2\le C',\ t\ge T,
\end{equation}
 where the constant $C'$ depends only on the constant $C$ in \eqref{4.smattr}. Moreover, this function solves the equation
 \begin{equation}\label{4.app}
 \eb\tilde u_{tt}+\tilde u_t+A(A\tilde u+f(\tilde u))=g+\varphi(t),
 \end{equation}
 with
 \begin{equation}\label{4.phi}
 \|\varphi\|_{C_b([T,+\infty),H)}\le C\sigma^\kappa,
 \end{equation}
 where $C$ and $\kappa>0$ are independent of $\sigma$
 and such that
 \begin{equation}\label{4.close}
 \|u(t)-\tilde u(t)\|_H+\|u_t(t)-\tilde u_t(t)\|_{H^{-2}}\le \sigma,
 \end{equation}
 for all $t\ge T$.
 \end{lemm}
 \begin{proof} The proof of this lemma is very similar to
 the proof of Theorem \ref{Zelik}. However, instead of the
 backward attraction to the smooth set of equilibria, we need to use the forward attraction to the smooth trajectory
 attractor $\mathcal A^{tr}_\eb$. Indeed, since the trajectory $(u,u_t)$ is attracted by the trajectory attractor $\mathcal A_\eb^{tr}$,
 then, analogously to \eqref{conv}, we have that, for every $\sigma>0$,
 there exist $T=T(\eb,\sigma,M_u(0))$ and a trajectory
 $(u_s,\Dt u_s)\in\mathcal K_\eb$, $s\ge1$ such that
 \begin{equation}\label{4.aa}
 \|u(t)-u_s(t)\|_{H}+\|\Dt u(t)-\Dt u_s(t)\|_{H^{-2}}\le \sigma,\ \ t\in[T+s-1,T+s+2]
 \end{equation}
 for every $s\ge1$ (compare with \eqref{zelik3}). Moreover, analogously to Theorem \ref{Zelik},
  due to the regularity \eqref{4.smattr} of the attractor and estimate \eqref{4.aa}, we  conclude that
\begin{equation}\label{4.estt}
 \|\Dt u_s(t)-\Dt u_{s+1}(t)\|_{H^1}+\|u_s(t)-u_{s+1}(t)\|_{H^3}\le C\sigma^\kappa,\ \ t\in[T+s,T+s+1],
\end{equation}
 where the positive constants $C$ and $\kappa$ are independent of $\eb$, $s$ and the choice of the trajectory $u$. Finally,
 defining the function $\tilde u(t)$ on the interval $[T+N,T+N+1]$ as follows
 $$
 \tilde u(t):=\theta(t-T-N)u_{N+1}(t)+(1-\theta(t-T-N))u_{N}(t),
 $$
 where the cut-off function $\theta$ is the same as in Theorem \ref{Zelik} and setting
 $$
 \varphi (t):=\eb\tilde u_{tt}(t)+\tilde u_t(t)+A(A\tilde u(t)+f(\tilde u(t)))-g,
 $$
 one can  see that the function $\varphi(t)$ satisfies estimate \eqref{4.phi}. Indeed, since $u_N$ and $u_{N+1}$ solve the initial problem \eqref{CH1-regular}, we have
 \begin{align}\label{4.tail}
 \varphi(t)&=2\eb\theta^\prime\Dt(u_{N+1}-u_N)
 +(\eb\theta''+\theta')(u_{N+1}-u_N)\\
 \nonumber
 &+A[f(\theta u_{N+1}+(1-\theta)u_{N})-\theta f(u_{N+1})-(1-\theta)f(u_{N})],
 \end{align}
 where $\theta=\theta(t-T-N)$. We see that the first two terms are immediately under the control thanks to \eqref{4.estt} and, in order to estimate the third term, we transform it as follows:
 \begin{align*}
 &f(\theta u_{N+1}+(1-\theta)u_{N})-\theta f(u_{N+1})-(1-\theta)f(u_{N})\\
 &=\theta(1-\theta)(u_{N+1}-u_{N})
 \int_0^1[f'(s\tilde u+(1-s)u_{N})-f'(s\tilde u+(1-s)u_{N+1})]\,ds.
 \end{align*}
 Thus, due to \eqref{4.estt} and the fact that $f\in C^4(\RR;\RR)$,
 the third term in \eqref{4.tail} is also under control and estimate \eqref{4.phi} holds.
  That finishes the proof of the lemma.
 \end{proof}

The next theorem is analogous to Theorem \ref{Zelik1}, but a bit more delicate since the regularity \eqref{4.phi}
of the function $\varphi$ given in Lemma \ref{Lem4.appr} is slightly lower than in Theorem \ref{Zelik}.

\begin{teor}\label{Th4.main} Let the assumptions of Lemma \ref{Lem4.appr} hold.
Then, for any $L>0$ and any trajectory $U\in K^+_\eb$, there exists
a time $T=T(\eb, L,M_u(0))$ such that the equation
\begin{equation}\label{4.smooth}
\eb v_{tt}+v_t+A(Av+f(v))+LA^{-1}v=G(t):=g+LA^{-1}u(t)
\end{equation}
has a regular  global solution $v(t)$, $t\ge T$ satisfying
\begin{equation}\label{4.Lreg}
\eb\|v_{tt}(t)\|^2_H+\|v_t(t)\|_{H^2}^2+\|v(t)\|_{H^4}^2\le C_L,
\end{equation}
where the constant $C_L$ depends on $L$, but is independent of $\eb$, $t$ and $u$. Moreover,
\begin{equation}\label{4.basic}
\|v_t(t)\|_{L^\infty(\Omega)}^2+\|v(t)\|_{L^\infty(\Omega)}^2\le C,
\end{equation}
where $C$ is independent of $L$ and $\eb$.
\end{teor}
\begin{proof} Step 1. As in Theorem \ref{Zelik1}, in order to solve equation \eqref{4.smooth} we introduce the function $w(t):=v(t)-\tilde u(t)$, where $\tilde u(t)$ is constructed in Lemma \ref{Lem4.appr}. Then, this function satisfies
\begin{equation}\label{4.eqw}
\eb w_{tt}+w_t+A(Aw+[f(\tilde u(t)+w)-f(\tilde u(t))])+LA^{-1}w=\tilde G(t):=LA^{-1}(u(t)-\tilde u(t))-\varphi(t).
\end{equation}
Moreover, due to Lemma \ref{Lem4.appr}, for any $L$ and any $0<\sigma<1/L^2$, we may find $T=T(\eb,L,\sigma,M_u(0))$ such that
\begin{equation}\label{4.dist}
\|\tilde G\|_{C_b([T,+\infty),H)}\le C\sigma^\kappa,
\end{equation}
where the positive constants $C$ and $\kappa$ are independent of $L$, $t$ and $u$.
Then, for small $\sigma$, equation \eqref{4.eqw} endowed by the initial data
\begin{equation}
\label{inidataw}
w(T)=0,\ \ w_t(T)=LA^{-1}(u(T)-\tilde u(T)),
\end{equation}
can be uniquely solved using the inverse function theorem (analogously to Theorem \ref{Zelik1}). However,
since we now do not have the control of the time derivative of $\varphi$, we are unable to construct the $\calV_3$-solutions
in such way and should restrict ourselves to consider the $\calV_1$-solutions only. Namely, it is not difficult to show
using the inverse function theorem that, for sufficiently small $\sigma$, there exists a unique global solution $w(t)$
of problem \eqref{4.eqw} with the above initial data such that
\begin{equation}\label{4.dist1}
\eb\|w_t(t)\|_{H}^2+\|w(t)\|^2_{H^2}\le C'\sigma^{2\kappa},\ \ t\ge T,
\end{equation}
where the constant $C'$ is independent of $L$, $t$ and $\sigma$
(since very similar arguments have been considered in detail in the
proof of Theorem \ref{Zelik1}, we omit the details here).
\par
In addition, from equation \eqref{4.eqw}, we see that
\begin{equation} \label{4.ode}
\eb A^{-1}w_{tt}+A^{-1}w_t=h_w(t):=
-Aw(t)-[f(\tilde u(t)+w(t))-f(\tilde u(t))]+A^{-1}\tilde G(t).
\end{equation}
Thus, due to estimates \eqref{4.dist} and \eqref{4.dist1}, we have
$$
\|h_w(t)\|_{H}\le C\sigma^\kappa,\quad t\ge T.
$$
Solving explicitly \eqref{4.ode} as an ODE with respect to $w_t$ and
using the last estimate together with \eqref{4.close} and
\eqref{inidataw}, we arrive at
\begin{equation}\label{4.dist2}
\|w_t(t)\|_{H^{-2}}\le C\sigma^\kappa,\quad t\ge T,
\end{equation}
where the positive constants $C$ and $\kappa$ are independent of $L$, $\eb$, $t$ and $u$.
\par
Step 2. As we have already mentioned, the regularity \eqref{4.dist1}
and \eqref{4.dist2} for the auxiliary problem \eqref{4.smooth} is
not sufficient for our purposes and we need to improve it. To this end,
we remind  that, by the construction of $\tilde u$, we may assume
without loss of generality that $\varphi (T)=0$ and, consequently,
initial conditions \eqref{inidataw} imply that $w_{tt}(T)=0$.
Therefore,
\begin{equation}\label{4.init}
v(T)=\tilde u(T),\quad v_t(T)=\tilde u_t(T)+LA^{-1}(u(T)-\tilde u(T)),
\quad v_{tt}(T)=\tilde u_{tt}(T).
\end{equation}
Thus, the initial data for the solution $v$ at $t=T$ is more regular,
namely (cf. \eqref{4.sm} and \eqref{4.close}),
\begin{equation}\label{4.iniest}
\|v(T)\|_{H^4}^2+\|v_t(T)\|_{H^2}^2+\eb\|v_{tt}(T)\|_{H}^2\le C,
\end{equation}
where the constant $C$ is independent of $\eb$ and $L$ (recall that $\sigma\in (0,1/L^2)$). Since the
(global in time and independent of $\eb$) control of the $H^2$-norm
of the solution $v(t)$ is already obtained (see estimates
\eqref{4.sm} and \eqref{4.dist1}), then $f(v(t))$ is under control.
Thus the regularity \eqref{4.iniest} of the initial data implies in a
standard way (see e.g., \cite{GSSZ} and the proof of Theorem
\ref{Zelik1}), that the solution $v(t)$ is indeed more regular and, in
particular, estimate \eqref{4.Lreg} holds with the constant $C_L$
depending on $L$, but independent of $\eb$ and $t$.
\par
Step 3. We only need to obtain the estimate \eqref{4.basic} for the
$L^\infty$-norms of $v$ and $v_t$ with the constant $C$ {\it
independent} of $L$. Moreover, the desired estimate for $v$ is
already available (see \eqref{4.dist1} and \eqref{4.sm}), so we
only need to estimate the $L^\infty$-norm of the time derivative.
To this end, we use the following standard interpolation inequality
together with \eqref{4.sm} and \eqref{4.dist2}
\begin{align}\label{4.mult}
&\|v_t(t)\|_{L^\infty}\le \|\tilde u_t(t)\|_{L^\infty}+\|w_t(t)\|_{L^\infty}\\
\nonumber
&\le C+\|w_t(t)\|_{H^{-2}}^{1/8}\|w_t(t)\|_{H^2}^{7/8}\le C+C\sigma^{\kappa/8}(1+\|v_t(t)\|_{H^2}+\|\tilde u_t(t)\|_{H^2})\le
C+C_L\sigma^{\kappa/8},
\end{align}
where only the constant $C_L$ may depend on $L$ and all of the constants are independent of $\eb$ and $\sigma$. It only remains
to note that, for every fixed $L$, we may fix $\sigma=\sigma(L)$ in a such way that $C_L\sigma^{\kappa/8}\le1$. Estimate \eqref{4.basic} is then an immediate corollary of \eqref{4.mult} and the theorem is proved.
\end{proof}
We are now ready to state and prove the main result of this section.

\begin{teor}\label{Th4.main1} Let the assumptions of Lemma \ref{Lem4.appr} hold, so that
the trajectory attractor $\mathcal A_\eb^{tr}$ of equation
\eqref{CHin}  satisfies the regularity estimate \eqref{4.smattr} with
the constant $C$ independent of $\eb\to0$.  Then, there exist a
$R_0$-ball
\begin{equation}\label{4.ball}
\mathcal B:=\{(a,b)\in H^4\times H^2,\ \ \|a\|^2_{H^4}+\|b\|_{H^2}^2\le R_0^2\},
\end{equation}
where the radius $R_0$ is independent of $\eb\to0$, and a
monotone function $Q_\eb$ (depending on $\eb$) such that, for
every weak energy solution $U(t)=(u(t),u_t(t))\in K_\eb^+$, the
following estimate holds:
\begin{equation}\label{4.expattr}
\dist_{\calX_0}(U(t),\mathcal B)\le Q_\eb(M^\eb_u(0))e^{-\beta t},
\end{equation}
where the positive exponent $\beta$ is independent of $t$ and $\eb\to0$.
\end{teor}

\begin{proof} Let $U\in K^+_\eb$ be arbitrary. Then, due to Theorem \ref{Th4.main}, for every $L>0$, there exists
$T=T(\eb,L,M^\eb_u(0))$ (which is independent of the concrete choice of the trajectory $U$) such that the auxiliary
equation \eqref{4.smooth} is  globally solvable for $t\ge T$ in the class of regular solutions and the solution $v(t)$
satisfies estimates \eqref{4.Lreg} and \eqref{4.basic}. Moreover, crucial for our method is the fact that the constant $C$ in
\eqref{4.basic} is independent  of $L$, $\eb$, $T$ and $U$. Let now $w(t):=u(t)-v(t)$. Then, this function solves
\begin{equation}\label{4.dif}
\eb w_{tt}+w_t-A(Aw+[f(v+w)-f(v)])+LA^{-1}w=0.
\end{equation}
Recall that we deal only with energy solutions which {\it can be
obtained by Galerkin approximations}. So, we now need to derive the
proper estimate for the function $w$ using the Galerkin
approximations exactly as in Theorem \ref{backsmooth}. However,
in order to avoid the technicalities, we will proceed by formal
multiplication of the equation \eqref{4.dif} by $A^{-1}(w_t+\alpha
w)$ (see the proof of Theorem \ref{backsmooth}  for the
justification). Then, after integration in space, we arrive at
\begin{equation}
\label{4.energy3}
\frac{d}{dt} E_{w} + \alpha E_{w} = H_{w}, \qquad\textrm{ in } [T,+\infty),
\end{equation}
where
\begin{align}
\label{4.energyk}
E_{w} &=  \eps\Vert A^{-1/2}w_t\Vert_H^2 + \Vert A^{1/2} w\Vert_H^2
+ L\Vert A^{-1} w\Vert_H^2 + 2\alpha\epsi( w, w_t) \\
\nonumber
&+ \alpha\Vert A^{-1/2}w\Vert_H^2
+ 2\left(F(v + w) - F(v) - f(v)w,1\right)
\end{align}
and

\begin{align}
\label{4.sourcek}
H_{w} &:=  -(2-3\eps)\Vert A^{-1/2} w_t\Vert_H^2
-\alpha \Vert A^{1/2} w\Vert_H^2
- \alpha L\Vert A^{-1} w\Vert_H^2\\
\nonumber
&+2\alpha\left(F(v + w) - F(v) - f(v)w
- (f(v+ w) - f(v))w,1\right)\\
\nonumber
&+ 2\alpha^2\epsi( w, w_t) + \alpha^2\Vert w\Vert_H^2
\\
\nonumber
&+ 2\left(f(v + w) - f(v) - f'(v)w,v_t\right).
\end{align}
Since the constant in estimate \eqref{4.basic}  for the
$L^\infty$-norm of $v$ and $v_t$ is independent of $L$, arguing as
in estimates \eqref{sourcekest}, we may fix $L$ in such way that (recall also \cite[(2.51)]{Ze04})
\begin{equation}\label{4.goodest}
 C \ge E_w\ge \gamma\|(w,w_t)\|_{\calX_0}^2,\qquad  H_w\le0,
\end{equation}
where the positive constants $C$ and $\gamma$ are independent of the concrete choice of the trajectory $U$.
\par
The Gronwall inequality applied to \eqref{4.energy3} now gives
$$
\|(w(t),w_t(t))\|_{\calX_0}^2\le C([M^\eb_u(T)]^2+1)e^{-\alpha (t-T)},\ \ t\ge T.
$$
It only remains to recall that $T=T(\eb,M_u^\eb(0))$ (and $L$ is now fixed) and use \eqref{M2} in order to
estimate $M_u^\eb(T)$ through $M^\eb_u(0)$. That yields
$$
\|(w(t),w_t(t))\|_{\calX_0}^2\le Q_\eb(M_u^\eb(0))e^{-\alpha t},
$$
which implies \eqref{4.expattr} and finishes the proof of the theorem.
\end{proof}

\begin{osse}\label{Rem4.bad} Although all
estimates of auxiliary solutions in the proof given above are uniform
with respect to $\eb\to0$, the monotone function $Q_\eb$ in the
main estimate \eqref{4.expattr} {\it depends} on $\eb$ since our
construction depends crucially on the attraction property to the
trajectory attractor $\mathcal A_\eb^{tr}$ in a weaker topology
(see Corollary \ref{attraction}) and the rate of that convergence may
be not uniform with respect to $\eb$. In a fact, it can be non-uniform
even with respect to $\eb\in[\eb_1,\eb_2]$ with $\eb_1>0$. This
problem may be solved in a standard way if we consider the
extended trajectory semigroup
$$
K_{[\eb_1,\eb_2]}^+:=\{(U,\eb),\ U\in K_\eb^+,\ \eb\in[\eb_1,\eb_2]\},\qquad
\mathbb T_\ell(U,\eb):=(\mathbb T_\ell U,\eb),
$$
construct its attractor and use the rate of  convergence to the
attractor for that semigroup. The difference is that we are now able
to approximate the trajectory $U\in K_{\eb}^+$ not only by the
elements of $\mathcal A_{\eb}^{tr}$ but also by the elements of
$\mathcal A_{\eb_n}^{tr}$ with $\eb_n\to\eb$, which would be
enough to obtain the uniformity with respect to $\eb$ (see \cite{BV}
for details). In order to avoid technicalities, we prefer not to give
the proof of that uniformity here. However, there is one more
essential drawback in the above scheme which cannot be corrected in
such an easy way, namely, as we have already mentioned, the
function $Q_\eb$ in estimate \eqref{4.expattr} depends in a crucial
way on the rate of attraction to the attractor and this rate of
convergence cannot be found explicitly or expressed in terms of the
physical parameters of the system considered (the usual drawback of
the global attractors theory), thus, we are factually unable to give
any expression for the function $Q_\eb$ following the above
described arguments. For this reason, we give below an alternative
explicit construction of the auxiliary ``almost solution'' $\tilde u(t)$
for the case when $\eb>0$ is very small. In addition, this
construction has another advantage, namely, it is based only on the
perturbation arguments and does not use the global Lyapunov
functional. This allows to apply it also to non-autonomous cases or to
the case of unbounded domains where the Lyapunov functional does
not exist. We will return to such issues in more details elsewhere.
\end{osse}
The following lemma is the uniform (with respect to $\eb\to0$) analogue of the key Lemma \ref{Lem4.appr}.
\begin{lemm}\label{Lem4.cappr}
Let the assumptions of Lemma \ref{Lem4.appr} hold. Then, for every
$\sigma>0$, there exists $\eb_0=\eb_0(\sigma)>0$ such that, for
any energy solution $U\in K^+_\eb$ with $\eb\le\eb_0$, there
exist $T=T(\kappa, M^\eb_u(0))$
 and a function $\tilde u(t)$ such that
 \begin{equation}\label{4.csm}
 \|\tilde u(t)\|_{H^4}^2+\|\tilde u_{t}(t)\|_{H^2}^2+\eb\|\tilde u_{tt}(t)\|_{H}^2\le C,\quad t\ge T,
\end{equation}
 where the constant $C$ is independent of $u$ and $\eb$ . Moreover, this function solves the equation
 \begin{equation}\label{4.capp}
 \eb\tilde u_{tt}+\tilde u_t+A(A\tilde u+f(\tilde u))=g+\varphi(t),
 \end{equation}
 with
 \begin{equation}\label{4.cphi}
 \|\varphi\|_{C_b([T,+\infty),H)}\le C\sigma^\kappa,
 \end{equation}
 where $C$ and $\kappa>0$ are independent of $\sigma$ and $\eb$
 and such that
 \begin{equation}\label{4.cclose}
 \|u(t)-\tilde u(t)\|_H+\|u_t(t)-\tilde u_t(t)\|_{H^{-2}}\le \sigma,
 \end{equation}
 for all $t\ge T$. In addition, all of the constants can be expressed explicitly in terms of the physical parameters.
 \end{lemm}
 \begin{proof} We define the trajectory $\tilde u$ as a solution of the modified Cahn-Hilliard equation
 \begin{equation}\label{4.limit}
 \tilde u_t+A(A\tilde u+f(\tilde u))+L_0A^{-1}\tilde u=g+L_0A^{-1}u(t),\ \ \tilde u(0)=u(0),
 \end{equation}
 where $u\in K^+_\eb$ and $L_0=L_0(f)$ is a sufficiently large number depending only on $f$. Then, on the one hand,
 the difference $u(t)-\tilde u(t)$ satisfies the following estimate:
 \begin{equation}\label{4.difpar}
 \|u(t)-\tilde u(t)\|_{H^{-1}}^2\le C\eb \(1+Q\(M^\eb_u(0)\)e^{-\alpha t}\),
 \end{equation}
 where the positive constants $C$ and $\alpha$ and the monotone function $Q$ are
 independent of $u$, $\eb$ and $t$
 (see \cite[Prop. 3.4]{GSSZ} for the details).
 Moreover, since we have the uniform control of the $H^1$-norms of $u$ and $\tilde u$, this estimate together with the interpolation inequality gives
 \begin{equation}\label{4.difpar1}
 \|u(t)-\tilde u(t)\|_{H}^2\le C\eb^{1/2}\(1+Q\(M^\eb_u(0)\)^2e^{-\alpha t}\).
 \end{equation}
 \par
 On the other hand, equation \eqref{4.limit} is a
 classical parabolic Cahn-Hilliard equation whose solutions
 possess the standard parabolic smoothing property.
 Their regularity is restricted only by the regularity of the nonlinearity
 $f$ and the external forces $g+L_0A^{-1}u(t)$. In our case,
 we have $f\in C^3(\R,\R)$ and
 $u\in W^{1,\infty}(\R_+, H^{-1}(\Omega))\cap
 L^\infty(\R_+, H^1(\Omega))$. It is then not difficult to verify that
 this regularity is enough to establish the following estimate for $\tilde u(t)$:
 \begin{equation}\label{4.regest}
 \|\tilde u(t)\|_{H^4}+\|\tilde u_t(t)\|_{H^2}
   +\sqrt{\eb}\|\tilde u_{tt}(t)\|_{L^2}\le C_*
   +Q(M^\eb_u(0))\frac{1+t^N}{t^N}e^{-\alpha t},
 \end{equation}
 where the positive constants $C_*$, $N$
 and $\alpha$ and the monotone function $Q$ are
 independent of $u$, $\eb$ and $t$.
 Since the derivation of this estimate is standard
 (although a bit technical), we leave it to the reader (see, e.g., \cite[Lemma~2.13]{EMZ}
 for a similar argument).

 Thus, the function $\tilde u(t)$ solves equation \eqref{4.capp} with
 $$
 \varphi(t):=\eb\tilde u_{tt}(t)+L_0A^{-1}(u(t)-\tilde u(t)),
 $$
 and therefore, thanks to \eqref{4.difpar1} and \eqref{4.regest},
 \begin{equation}\label{4.phifin}
 \|\varphi(t)\|_{H}\le C\eb^{1/4}\(1+Q(M^\eb_u(0))\frac{1+t^N}{t^N}e^{-\alpha t}\).
 \end{equation}
 Finally, we only need to estimate the difference
 $v(t):=u(t)-\tilde u(t)$ in the proper norm.
 To this end, taking the difference between \eqref{CH1-regular} and \eqref{4.capp}, we derive that
 \begin{equation}\label{4.ode1}
 \eb v_{tt}+v_t=H(t):=-A^2v(t)-A[f(u(t))-f(\tilde u(t))]-\varphi(t).
\end{equation}
Moreover, using assumption \eqref{f3} on the nonlinearity $f$ and
the fact that the $L^{p+4}$-norms of $u$ and $\tilde u$ are under
control (due to the energy estimate), we conclude from estimate
\eqref{4.difpar1} that
$$
\|f(u(t))-f(\tilde u(t))\|_{L^1}
\le(Q(M^\eb_u(0))e^{-\alpha t}+C_*)\|u(t)-\tilde u(t)\|_H^\kappa\le C\eb^{\kappa/4}\( 1+Q(M^\eb_u(0))e^{-\alpha t}\)
$$
for some positive $\kappa$ depending only on $p$. Thus, thanks to
\eqref{4.difpar1} and \eqref{4.phifin} and the embedding
$H^2(\Omega)\hookrightarrow C^0(\bar\Omega)$,
$$
\|A^{-2}H(t)\|_H\le C\eb^{\kappa/4}\(1+Q(M^\eb_u(0))\frac{1+t^N}{t^N}e^{-\alpha t}\).
$$
Solving explicitly the ODE \eqref{4.ode1} on the interval $[t/2,t]$ and
using the last estimate together with the estimate \eqref{4.regest},
we conclude that
$$
\|A^{-2}v_t(t)\|_{H}\le C\eb^{\kappa/4}
+ C(\eb^{\kappa/4}+e^{-\frac t{2\eb}})\(1+Q(M^\eb_u(0))
\frac{1+t^N}{t^N}e^{-\alpha t}\).
$$
Finally, keeping in mind that the $H^{-1}$-norm of $v_t$ is under
control (due to the energy estimate), we end up with
\begin{equation}\label{4.lastu}
\|u_t(t)-\tilde u_t(t)\|_{H^{-2}}\le C\eb^{\kappa/12}+
C(\eb^{\kappa/12}+e^{-\frac t{6\eb}})\(1+Q(M^\eb_u(0))
\frac{1+t^N}{t^N}e^{-\alpha t}\).
\end{equation}
Estimates \eqref{4.difpar1}, \eqref{4.regest}, \eqref{4.phifin} and \eqref{4.lastu} show that, indeed, for every $\kappa>0$, we may fix $\eb_0=\eb_0(\kappa)$ and $T=T(M^\eb_u(0))$ such that all estimates stated in Lemma \ref{Lem4.cappr} will be satisfied uniformly with respect to $\eb\le \eb_0$ and $U\in K^+_\eb$. That finishes the proof of the lemma.
 \end{proof}
\begin{coro}\label{Cor4.e0} Let the assumptions of Lemma \ref{Lem4.appr} hold. Then, there exist $\eb_0>0$ and $R_0>0$ such that, for every $\eb<\eb_0$, the set
\begin{equation}\label{4.uniball}
\mathcal B_\eb:=\{(u_0,u_0')\in H^4\times H^2: \|u_0\|^2_{H^4}+\|u_0'\|_{H^2}^2 \le R_0^2\}
\end{equation}
attracts exponentially all energy solutions $U\in K_\eb^+$:
\begin{equation}\label{4.uniest}
\dist_{\calX_0}\((u(t),u_t(t)),\mathcal B_\eb\)\le Q(M^\eb_u(0))e^{-\alpha t},
\end{equation}
where the positive constant $\alpha$ and monotone function $Q$ are
independent of $t$, $\eb\le\eb_0$ and of the concrete choice of the
trajectory $U\in K_+^\eb$.
\end{coro}
Indeed, the derivation of estimate \eqref{4.uniest} is analogous to
the proof of Theorem \ref{Th4.main1} with the only difference that,
instead of the non-uniform approximations of Lemma
\ref{Lem4.appr}, one should use the uniform approximations of
Lemma \ref{Lem4.cappr}.

\begin{osse}\label{Rem4.M} In fact, we have proven a
bit more than  \eqref{4.expattr} or \eqref{4.uniest}. Namely,
recalling \cite[Def.~4.1]{Ze04} (see also \eqref{Mfunctional}), let
us consider the $M$-distance to the set $\mathcal B_\eb$ defined
by
\begin{align}
&\dist_{M_u^\eb}(t,\mathcal B_\eb) \\
\nonumber
&:=\inf \Big\{\liminf_{k\to
\infty} d_{\calX_0}(U^{n_k}(t),P_{n_k}\mathcal B_\eb) \, : \, U
= \Theta^+ - \lim_{k\to \infty}U^{n_k}, \; U(0) = [\calX_0]^w  -
\lim_{k\to \infty}U^{n_k}(0) \Big\}.
\end{align}
Recall that the external infimum is taken over all the sequences $\{U^{n_k}(t)\}_{k\in {\mathbb N}}$ of
Faedo-Galerkin approximations which
$\Theta^+$-converge to the given solution $U$. Then, we may improve estimate \eqref{4.uniest} as follows:
\begin{equation}\label{4.Mest}
\dist_{M_u^\eb}(t,\mathcal B_\eb)\le Q(M^\eb_u(0))e^{-\alpha t}.
\end{equation}
This slight generalization is however important for applying the arguments based on the transitivity of exponential attraction, see \cite{FGMZ} and next section.
\end{osse}
\section{Exponential attractors for energy solutions}\label{s.exp}
In this concluding section, we discuss the exponential attractors for
problem \eqref{CH1-regular}. We start with the case of small $\eb$
(and the 3D case for definiteness). We first recall that, due to
Theorem \ref{globexuni} (see \cite[Thm.~2.7]{GSSZ}), equations
\eqref{CH1-regular} are globally solvable in the class of more
regular $\calV_1$-solutions if the initial data is not large enough and
$\eb$ is small enough. More precisely, the equation generates a
dissipative semigroup $S_\eb(t)$ on the set $\mathbb
B_\eb\subset \calV_1$ defined by \eqref{inclB} and, in particular,
the phase space $\mathbb B_\eb$ of that semigroup contains an
$R(\eb)$-ball of the space $\calV_1$ with $R(\eb)\to\infty$ as
$\eb\to0$.
\par

\begin{osse}\label{robust}
It is worth observing that, on account of the results obtained in \cite{GSSZ},
we can argue as in \cite{FGMZ} to prove that
the family  $\{(\mathbb B_\eb,S_\eb(t))\}_{\eb\in [0,\eb^\prime]}$, for some $\eb^\prime>0$ small enough,
possesses a uniform family of
exponential attractors $\mathcal M_\eb\subset \calV_3 \cap (H^4\times H)$
with the following properties:
\par
\item{(i)} The sets $\mathcal M_\eb$ are uniformly bounded in $H^4\times H$ as $\eb\to0$.
\par
\item{(ii)} The sets $\mathcal M_\eb$  are compact in $H^3\times H^{-1}$  and their fractal dimensions are uniformly bounded, i.e.,
$$
\dim_f(\mathcal M_\eb, H^3\times H^{-1})\le C,
$$
where $C$ is independent of $\eb\to0$.
\par
\item{(iii)} The uniform exponential attraction property holds
\begin{equation}\label{5.regexp}
\dist_{\calV_1}(S_\eb(t)\mathbb B_\eb,\mathcal M_\eb)\le Ce^{-\alpha t}
\end{equation}
with positive $C$ and $\alpha$ independent of $\eb\to0$.
\par
\item{(iv)} $\mathcal M_\eb$ tends to the limit exponential attractor $\mathcal M_0$ as $\eb\to0$ in the following sense
$$
\dist^{symm}_{H^3\times H^{-1}}\(\mathcal M_\eb,\mathcal M_0\)\le C\eb^\kappa
$$
for some positive $\kappa$ and $C$ which are independent of $\eb$.
We also recall that, as usual, in order to compare the solutions of the hyperbolic equation \eqref{CH1-regular} with $\eb>0$ and the solutions of the limit parabolic problem which corresponds to $\eb=0$, one needs to extend the limit parabolic semigroup to the surface (see \cite{BV, FGMZ}  for details)
$$
\mathbb S:=\{(u,v)\in H^4\times H,\ \ v=-A(Au+f(u))+g\}.
$$
\end{osse}

\par
The aim of this section is to verify that the above exponential
attractors $\mathcal M_\eb$ attract exponentially not only the weak
solutions (see \eqref{5.regexp}), but also all energy solutions of
problem \eqref{CH1-regular}. Namely, the following theorem holds.

\begin{teor}\label{Th5.expattr} Let the assumptions of Theorem \ref{existener} hold. Then, there exists $\eb_0>0$ such that,
for every $\eb<\eb_0$, there exists a family of exponential attractors $\mathcal M_\eb$ satisfying the properties 1)-4) formulated above and, in addition, for every energy solution $U\in K_\eb^+$ of problem \eqref{CH1-regular},
\begin{equation}\label{5.enerexp}
\dist_{\calX_0}\((u,u_t),\mathcal M_\eb\)\le Q(M^\eb_u(0))e^{-\alpha t},
\end{equation}
where the positive constant $\alpha$ and monotone function $Q$ are independent of $\eb$, $t$ and $u$.
\end{teor}
\begin{proof} Indeed, according to Corollary \ref{Cor4.e0} and Remark \ref{Rem4.M},
\begin{equation}\label{5.1}
\dist_{M^\eb_u}(t,\mathcal B_\eb) \le Q(M^\eb_u(0))e^{-\alpha t},
\end{equation}
where $Q$ and $\alpha$ are independent of $\eb\le\eb_0$, $t$ and $u$. On the other hand, thanks to \eqref{5.regexp},
\begin{equation}\label{5.2}
\dist_{\calV_1}(S_\eb(t)\mathcal B_\eb,\mathcal M_\eb)\le Ce^{-\alpha t}
\end{equation}
if $\eb>0$ is small enough (so that $\mathcal B_\eb \subseteq \mathbb{B}_{\eps}$, cf. \eqref{inclB}).
Then, keeping in mind  that
$\calV_1\subset\calX_0$, we see that, in order to prove estimate
\eqref{5.enerexp} through the transitivity of the exponential attraction
(see \cite[Thm.~5.1]{FGMZ}, cf. also \cite[Sec.~4]{Ze04}), we
only need to check the following version of Lipschitz continuity (see
\cite[(5.1)]{FGMZ}):
\begin{equation}\label{5.3}
\dist_{M^\eb_u}(t,V(t))\le \dist_{M^\eb_u}(0,V(0))e^{Kt},
\end{equation}
where $\eb>0$ is small enough, $U\in K_\eb^+$, $V(t):=S_\eb(t)V_0$, $V_0\in\mathcal B_\eb$, is an arbitrary  strong solution of equation \eqref{CH1-regular} starting from the set $\mathcal B_\eb$ and the positive constant $K$ is independent of $\eb$, $V$ and $t$.
This Lipschitz continuity property can be easily verified arguing as in the proof of Theorem \ref{unismooth}. Indeed, since (due to Theorem \ref{globexuni}) the solution $V(t)=(v(t),v_t(t))$ exists globally and satisfies the dissipative estimate, one can show that
$$
\|v(t)\|_{H^4}^2+\|v_t(t)\|^2_{H^2}+\eb\|v_{tt}(t)\|_{H}^2\le R_1=Q(R_0),
$$
where the constant $R_1$ is independent of $\eb$, $V(0)\in\mathcal B_\eb$ and $t$. Therefore, we have the uniform control
\begin{equation}\label{5.c}
\|v(t)\|_{L^\infty(\Omega)}+\|v_t(t)\|_{L^\infty(\Omega)}\le C,
\end{equation}
where $C$ is independent of $\eb$, $t$ and $V(0)$. Thus, defining $V^{n_k}(t)=P_{n_k}V(t)$, $z^{n_k}:=u^{n_k}-v^{n_k}$, where
$u^{n_k}(t)$ are the Faedo-Galerkin approximations to the solution $U\in K^+_\eb$ (see Section \ref{sectraj}) and arguing exactly as in the proof of Theorem \ref{unismooth}, we derive that
\begin{equation}
\label{5.en}
\frac{d}{dt} E_{z^{n_k}}(t) + \alpha E_{z^{n_k}}(t) \le C\|h^{n_k}(t)\|_{H^{-1}}^2,
\end{equation}
where $E_{z^{n_k}}$ and $h^{n_k}$ are defined as in \eqref{energyk} and \eqref{hhhh} (with $\tilde u$ replaced by $v$) respectively. On the other hand, we have
$$
\|h^{n_k}(t)\|_{H^{-1}}^2\le C\|P_{n_k}(f(v(t))-f(v^{n_k}(t))\|_{H^1}^2+CL^2E_{z^{n_k}}(t),
$$
and applying the Gronwall inequality to \eqref{5.en}, we have
$$
E_{z^{n_k}}(t)\le E_{z^{n_k}}(0)e^{Kt}
+C\int_0^te^{K(t-s)}\|P_{n_k}(f(v(s))-f(v^{n_k}(s)))\|_{H^1}\,ds,
$$
for some positive $C$ and $K$ independent  of $u$, $v$, $n_k$,
$\eb$ and $t$. Passing now to the limit $k\to\infty$ and using that
the right-hand side tends to zero (note that $v$ is smooth), we derive the
desired Lipschitz continuity \eqref{5.3}. Estimate \eqref{5.enerexp} is now a standard
corollary of transitivity of exponential attraction. This finishes the
proof of the theorem.
\end{proof}

We now consider the 2D case  with the growth restriction
\eqref{f3bis} for the nonlinearity $f$ (at most cubic growth rate).
Then, on account of \cite[Thms.~2.2, 3.1 and 5.1]{GSZ}), the 2D
problem \eqref{CH1-regular} generates a dissipative semigroup
$S_\eb(t)$ in the phase space $\calV_2$ for every finite $\eb>0$
and this semigroup possesses an exponential attractor $\mathcal
M_\eb$ which is bounded in $\calV_3$. Next theorem shows that
this exponential attractor attracts exponentially the energy solutions
as well.

\begin{teor}\label{Th5.2d} Let the assumptions of Theorem
\ref{smooth2D} hold. Then, for every $\eb>0$, the exponential
attractor $\mathcal M_\eb$ for the quasistrong $\calV_2$-solutions
constructed in \cite{GSZ} attracts exponentially energy solutions as
well. Namely, for any bounded set $B\subset\calX_0=\calV_0$, we
have
\begin{equation}\label{5.2dexp}
\dist_{\calX_0}(S_\eb(t)B,\mathcal M_\eb)
\le Q(\|B\|_{\calX_0})e^{-\alpha t},
 \end{equation}
 where the function $Q$ and constant $\alpha$
 are independent of $t$ and $B$, but may depend on $\eb$.
 Thus, $\mathcal M_\eb$ is an exponential attractor
 for the solution semigroup $S_\eb(t)$ acting on the
 energy phase space $\calX_0$ as well.
\end{teor}
Indeed, the proof of this theorem repeats word by word the proof of
the previous Theorem \ref{Th5.expattr}, with the only difference
that, instead of Corollary \ref{Cor4.e0}, one should use Theorem
\ref{Th4.main1}.
\par
To conclude, we apply the proved exponential regularization for the
2D case to one problem which remained unsolved in the previous
paper \cite{GSZ}. Namely, we have proved there that problem
\eqref{CH1-regular} with cubic growth restriction is well posed and
dissipative in $\calV_0$ and $\calV_2$, but the dissipativity in the
phase space $\calV_1$ occurred surprisingly more delicate and
remained an open issue. The next theorem fills this gap without any
restriction on $\eb$ (compare with \cite[Thm.~5.3]{CC}).

\begin{teor}\label{Th5.last} Let the assumptions of Theorem \ref{Th5.2d} hold. Then, for every $U(0)\in\calV_1$, problem \eqref{CH1-regular} is uniquely solvable in the phase space $\calV_1$ and the following dissipative estimate hold:
\begin{equation}\label{5.v1dis}
\|U(t)\|_{\calV_1}\le Q(\|U(0)\|_{\calV_1})e^{-\alpha t}+Q(\|g\|_{L^2}),
\end{equation}
where the positive constant $\alpha$ and monotone function $Q$ are independent of $t$ and $U(0)$.
\end{teor}
\begin{proof} Obviously, we only need to verify the dissipative estimate \eqref{5.v1dis}. As usual, we give only the formal derivation which can be easily justified using the Galerkin approximations. The first step is completely standard: we multiply equation \eqref{CH1-regular} by $u_t+\gamma u$ where $\gamma>0$ is a sufficiently small positive number and integrate over $x$. Then, after the straightforward transformations, we get
\begin{equation}\label{5.st}
\frac d{dt}Z_u(t)+\kappa Z_u(t)
+\kappa\|u_t\|^2_{L^2}\le C(1+\|g\|^2_{L^2}
+\|u(t)\|_{H^1}^2)+\frac{1}{2}|(f''(u)|\nabla u|^2,u_t)|,
\end{equation}
where $\kappa>0$ is small enough, $C>0$,
$$
Z_u(t):=\frac12\(\eb\|u_t\|^2_H+\|u\|^2_{H^2}
+(f'(u)\nabla u,\nabla u)+L_0\|u\|_H^2+2\gamma\eb(u,u_t)\),
$$
and the constant $L_0$ is chosen in such way that
\begin{equation}\label{5.eqpos}
k\|U(t)\|_{\calV_2}^2\le Z_u(t)\le Q(\|U(t)\|_{\calV_2})
\end{equation}
for some $k>0$ and some monotone function $Q$. Thus, the main
problem is how to estimate the last term in the right-hand side of
\eqref{5.st}.
\par
Actually, if we use the Br\'ezis-Gallouet logarithmic inequality together with the fact that
$|f''(u)|\le C(1+|u|)$ (cf.~\cite{BG},  see also
\cite[(2.34)]{GSZ}), we obtain
\begin{align}
\label{5.egg}
&|(f''(u)|\nabla u|^2,u_t)|\le \kappa\|u_t\|^2_H
+C(1+\|u\|_{L^\infty}^2)\|\nabla u\|^4_{L^4}
\\
\nonumber
&\le
C(1+\|u\|_{H^1}^2)\ln(1+\|u\|_{H^2}^2)\|\nabla
u\|^4_{L^4}+\kappa\|u_t\|_H^2\\
\nonumber
&\le
C'(1+\|u\|_{H^1}^2)\|\nabla
u\|^4_{L^4}\ln(1+Z_u(t))+\kappa\|u_t\|^2_H.
\end{align}
Notice now that, if we simply employ
$$
\|\nabla u\|^4_{L^4}\le C\|\nabla u\|^2_{H}
\|u\|^2_{H^2}\le C_1\|u\|^2_{H^1}Z_u(t),
$$
then (using also the $\calX_0$-energy estimate for the solution
$U$) we end up with an inequality of the form
$$
\frac d{dt}Z_u(t)+\kappa Z_u(t)\le CZ_u(t)\big(1+\ln(1+Z_u(t)\big)+C_1,
$$
which is enough to verify the global existence; however,  the
$\calV_2$-norm of the solution will diverge in time as a double
exponential. Thus, in order to obtain the desired dissipative estimate,
we have to proceed more carefully. Namely, thanks to Theorem
\ref{Th4.main1}, we can split the solution $u(t)=v(t)+w(t)$, where
\begin{equation}\label{5.split}
\|v(t)\|_{H^4}\le 2R_0,\ \ \ \|w(t)\|_{H^1}
\le Q(\|U(0)\|_{\calX_0})e^{-\alpha t},
\end{equation}
for some $R_0>0$, and estimate the $L^4$-norm of $\nabla u$ as
follows
\begin{align*}
&\|\nabla u\|^4_{L^4}\le 8(\|\nabla v\|^4_{L^4}+\|\nabla w\|^4_{L^4})\le C(R_0^4+\|w\|^2_{H^1}\|w\|^2_{H^2})
\\
&\le C(1+Q(\|U(0)\|_{\calX_0})e^{-\alpha t}\|u-v\|^2_{H^2})\le C+Q(\|U(0)\|_{\calX_0})e^{-\alpha t}(1+Z_u(t)).
\end{align*}
Inserting this estimate into the right-hand sides of \eqref{5.egg}
and \eqref{5.st} and using the dissipative estimate for the
$\calX_0$-energy norm of $U(t)$, we end up with the  refined
differential inequality
\begin{equation}\label{5.last}
\frac d{dt}Z_u(t)+\kappa Z_u(t)\le Q(\|U(0)\|_{\calX_0})e^{-\alpha t}Z_u(t)\ln(1+Z_u(t))+Q(\|g\|_H)+
Q(\|U(0)\|_{\calX_0})e^{-\alpha t}.
\end{equation}
It remains to note that (see proof of \cite[Thm.~3.1]{GSZ}) the
differential inequality \eqref{5.last} gives indeed the desired
dissipative estimate \eqref{5.v1dis} and finishes the proof of the
theorem.
\end{proof}
\begin{osse}\label{Rem5.last}
As we have already mentioned in Remark \ref{Rem4.bad}, the
functions $Q$ in Theorems \ref{Th5.2d} and \ref{Th5.last} depend
on the rate of convergence of weak energy solutions to the smooth
global attractor $\mathcal A_\eb$ and, by this reason, we cannot
find the explicit expressions of $Q$ in terms of the physical
parameters of the system. This drawback can be overcome if we construct the
smooth approximate solutions not by using the fact that the energy
trajectories tend to the smooth global attractor, but rather observing
that these trajectories visit regularly any arbitrarily small neighborhood
of the equilibria set (which follows from the existence of a global 
Lyapunov functional) and the smoothness of the set of equilibria.
However, this argument is much more delicate and, in order to avoid
the related technicalities, we will not present it here. 
\end{osse}




\vspace{35mm}

\noindent%
{\bf First author's address:}\\[1mm]
Maurizio Grasselli\\
Dipartimento di Matematica, Politecnico di Milano\\
Via E.~Bonardi, 9,~~I-20133 Milano,~~Italy\\
E-mail:~~{\tt maurizio.grasselli@polimi.it}

\vspace{4mm}

\noindent%
{\bf Second author's address:}\\[1mm]
Giulio Schimperna\\
Dipartimento di Matematica, Universit\`a degli Studi di Pavia\\
Via Ferrata, 1,~~I-27100 Pavia,~~Italy\\
E-mail:~~{\tt giusch04@unipv.it}

\vspace{4mm}

\noindent%
{\bf Third author's address:}\\[1mm]
Sergey Zelik\\
Department of Mathematics, University of Surrey\\
Guildford,~~GU2 7XH,~~United Kingdom\\
E-mail:~~{\tt S.Zelik@surrey.ac.uk}

\end{document}